\documentclass[%
  onecolumn
]{mpi2015-cscpreprint}

\usepackage[american]{babel}
\usepackage{tikz}
\usepackage{pgfplots}
\pgfplotsset{width=10cm,compat=1.9}
\usepackage[inkscape=true]{svg}
\usepackage{graphicx}
\usepackage[section]{placeins}
\usepackage{amssymb}
\usepackage{orcidlink}
\usepackage{fontawesome}
\usepackage{amsthm}
\usepackage{amsmath}
\usepackage{mathtools}
\usepackage{caption}
\usepackage{subcaption}
\usepackage{siunitx}

\usepackage{lineno}

\usepackage{filemod}
\usetikzlibrary{external}
\tikzexternaldisable% \includetikz locally activates it again, but we want it off, so it does not interfere with todonotes, etc.

\usepackage[normalem]{ulem}

\usepackage[textsize=scriptsize]{todonotes}
\usepackage{cleveref}
\usepackage[software, hardware]{mymacros}
\usepackage{siunitx}

\renewcommand{\hat}[1]{{\widehat{#1}}}
\renewcommand{\tilde}[1]{{\widetilde{#1}}}
\newcommand{\linew}{1pt}
\newcommand{\plotw}{2.4in}
\newcommand{\ploth}{2.3in}
\newcommand{\plotwtripple}{2.1in}
\newcommand{\plothtripple}{2.1in}
\begin{document}
  
%%%%%%%%%%%%%%%%%%%%%%%%%%%%%%%%%%%%%%%%%%%%%%%%%%%%%%%%%%%%%%%%%%%%%%%%%%%%%%%%
% PAPER INFORMATION.                                                           %
%%%%%%%%%%%%%%%%%%%%%%%%%%%%%%%%%%%%%%%%%%%%%%%%%%%%%%%%%%%%%%%%%%%%%%%%%%%%%%%%
\title{Reduced-Order Inference with Structure-Preserving Parametrization for Bending and Rotating Systems }
  
\author[$\dagger$,$\ast$, \faEnvelopeO]%
{Yevgeniya~Filanova~\orcidlink{0000-0002-8599-3747}}
  
\author[$\ast$]%
{\authorcr Igor~Pontes~Duff~\orcidlink{0000-0001-6433-6142}}

\author[$\ast$]%
{\authorcr Pawan~Goyal~\orcidlink{0000-0003-3072-7780}}

\author[$\ast$, $\dagger$]%
{Peter~Benner~\orcidlink{0000-0003-3362-4103}}
	
\affil[$\ast$]{%
		Max Planck Institute for Dynamics of Complex Technical Systems,
		Magdeburg, Germany}
\affil[$\dagger$]{Otto von Guericke University Magdeburg, Magdeburg, Germany}

\affil[\faEnvelopeO]{Corresponding author,\email{filanova@mpi-magdeburg.mpg.de}}


\abstract{%
Mechanical systems are often characterized only by their response to certain loads known from experiments or simulations. The obtained data can be used for various purposes: system analysis, design of mathematical models, or construction of reduced-order models for further simulations under different loading conditions. The use of data for reduced-order modeling is an important developing research direction, especially when the high-dimensional system operators are unknown and their low-dimensional approximation is required for accurate but fast simulations. Our goal is to obtain the low-dimensional surrogate model from the available input signal and deformation trajectory data, capturing the correct system behavior for the basic deformation cases, namely bending and rotation, which are present in almost every complex mechanical system.

In this work, we propose a methodology to infer the system operators by solving a nonlinear unconstrained optimization problem. The methodology is based on the operator inference approach for second-order systems and includes a parametrization of the unknown operators that preserves their symmetric positive definite or skew-symmetric structure. We demonstrate the performance of the novel approach for three numerical examples that are used to simulate basic bending and rotating.}

\novelty{\begin{itemize}
		\item Inference of a full physical representation of a mechanical system from simulation results and external input-signal data;
		\item  Formulation of an unconstrained nonlinear optimization problem with embedded structure-preserving parametrization of the unknown system matrices;
		\item Separate identification of the symmetric positive semi-definite damping and skew-symmetric gyroscopic terms.
	\end{itemize}}

\maketitle

\section{Introduction}%
\label{sec:intro}

% Mechanical systems: structural and rotation
Analysis and design of engineering structures are nowadays based on a large amount of different data from virtual and physical experiments. Empirical tests can be used, e.g., to develop mathematical models from input-output data relationship. Different techniques to determine models and predict the system output from the past observations are investigated in the system identification framework, as well as techniques to obtain coefficients, assuming a prior knowledge of the model structure \cite{bookIdentMech}. 
%For example, the state-space representation can be identified from time-domain data using the eigensystem realization algorithm (ERA) \cite{Jua1985}.

On the other hand, numerical experiments based on well-studied theoretical principles provide a large amount of data to design engineering processes. In case of mechanical systems, the simulations are often performed using the finite element method (FEM) \cite{Zienkiewicz2013}, implemented in well-established legacy codes. Accurate simulations of practical problems often require cumbersome models and fine meshing, which results in systems with large and sparse matrices. High computational costs limit the usage and efficiency of the simulation models, therefore constructing of surrogate models of much smaller size is required. Thus, we are interested in identification of small models from given data, which we call reduced-order inference.
%that can accurately represent the dynamics of the physical system encoded in 
% Non-intrusive MOR for second-order systems

Some of the system identification methods also provide small surrogate models from experimental measurements. For example, the eigensystem realization algorithm (ERA) utilizes impulse response data to construct linear discrete-time system \cite{era}. Another well-known technique called Loewner-framework provides the system matrices from frequency-domain data, i.e. from transfer function samples \cite{MAYO2007634}. The Loewner reduction approach can be used to obtain the state-space, as well as the second-order representation of a mechanical system \cite{PontesDuff2022}.

Often, it is necessary to identify the operators from time-domain data, including data sets that stem from numerical experiments. For these tasks, one of the important methods is the dynamic mode decomposition (DMD) \cite{morSch10}. The DMD approach operates with the time series of data and identifies the underlying dynamics of the given snapshots \cite{morKutBBetal16}. Reduction of the state is usually done via the singular value decomposition (SVD) of the snapshot matrix. There exist many extensions of the DMD method for linear and nonlinear dynamical systems, such as \cite{dmd_article, Benner2018, dmd_bai, morProBK16}.

% Existing non-intrusive methods and their drawbacks
Another method to identify the system operators from state snapshots is operator inference (OpInf) \cite{morPehW16}. It provides the continuous time operators via solving a least-squares problem for the low-dimensional approximation of the snapshot data. The inference of the reduced model is done in a form inspired by intrusive projection and preserves the structure of the original discretized-in-space governing system of equations. A review of the basic theoretical and practical aspects of operator inference for first-order differential equation systems can be found in \cite{Kramer2024}. There exist many extensions and adaptations of the OpInf method, e.g., for Navier-Stokes equations \cite{BenGHP22}, and second-order differential equation systems that usually model the mechanical system behavior. In \cite{SHARMA2024134128, SHARMA2024116865} the inference is based on the Lagrangian energy function and the system operators, with some restrictions for the mass matrix, are provided from state and input signal snapshots. In \cite{filanova2022operator} the inference is based on the discretized-in-space form of the second-order system of equations and provides the system matrices using state snapshots and available external force snapshots. One of the important aspects is preserving the stability of the inferred system. The regularization terms can be used to induce the stability bias to the system, which was done for quadratic models in  \cite{SAWANT2023115836}. Another way is introducing a constrained optimization problem to preserve the symmetric positive definite (SPD) structure of the matrices, as proposed in \cite{SHARMA2024134128, SHARMA2024116865,filanova2022operator}. 

However, obtaining all the second-order system operators without additional assumptions for the mass matrix and without using the external load snapshots is not fully resolved in the above mentioned methods. In this work, we propose a methodology that allows to infer separately the system matrices using only the input-signal data and the data available from numerical integration. Moreover, depending on the simulation case, we guarantee to infer the SPD mass, stiffness, damping matrices and the skew-symmetric gyroscopic matrix without constrained optimization. Following the ideas proposed in \cite{goyal2023stab, goyal2024quadstab} for continuous linear and quadratic systems, we introduce a parametrization of the unknown matrix coefficients, such that the SPD and skew-symmetric properties are guaranteed by construction. This leads to a nonlinear optimization problem that can be efficiently solved using automatic differentiation tools, providing very small reduced-order models.

%\cite{TisoMahdiabadi}
% parameterization of stable matrices \cite{GILLIS2017113}

% Section organization
The paper is organized as following. \Cref{sec:2order} contains a description of the systems that are suitable for the application of our method, which is described in \Cref{sec:popinf}. The latter section contains also practical details about the implementation of the presented approach. In \Cref{sec:numres} three numerical examples for bending and rotation structures are presented.

\section{Second-order systems} \label{sec:2order}

A mechanical system converts an input motion and force to the output motion and force. In this work, we consider linear mechanical systems, discretized in space, in the following form:

\begin{equation}
\label{eq:original_noG}
\textbf{M}\ddot{\mathbf{x}}(t) + \textbf{E}\dot{\mathbf{x}}(t) + \textbf{K}\mathbf{x}(t) = \textbf{B} \mathbf{u}(t),
\end{equation}
where $\textbf{M}$,$\textbf{K} \in \mathbb{R}^{n \times n}$ are the SPD mass and stiffness matrices, respectively, $\textbf{B}$ represents the input matrix, $\textbf{x}(t) \in \mathbb{R}^{n}$ stands for the displacement vector, $\textbf{u}(t) \in \mathbb{R}^{m}$ is the input vector. 

The matrix $\textbf{E} \in \mathbb{R}^{n \times n}$ is typically equal to the internal positive semi-definite damping matrix $\textbf{D}$ that approximates the viscous damping in the system. When the structure is rotating, the matrix $\textbf{E}$ also contains a skew-symmetric gyroscopic matrix $\textbf{G}$, which can include also Coriolis damping effects \cite{rotorG, Lal97, Baruh99}. In this case, the system of governing equations takes the following form:

\begin{equation}
\label{eq:original}
\textbf{M}\ddot{\mathbf{x}}(t) + \textbf{D}\dot{\mathbf{x}}(t) + \Omega \textbf{G} \dot{\mathbf{x}}(t) + \textbf{K}\mathbf{x}(t) = \textbf{B} \mathbf{u}(t),
\end{equation}
where $\Omega$ is the constant spin speed. Note that the system \eqref{eq:original} is equivalent to the system \eqref{eq:original_noG} with $\textbf{E} = \textbf{D} + \Omega \textbf{G}$.

As mentioned above, commercial software packages usually do not provide an access to the dynamical equations. Moreover, if the input matrix $\textbf{B}$ contains software-specific node-distribution information, etc., it can not be easily extracted from the user-defined external force information. Therefore, relying on the known system structure \eqref{eq:original} or \eqref{eq:original_noG}, we aim to identify a low-dimensional surrogate model, i.e., infer the reduced matrix coefficients $\widetilde{\textbf{M}},\widetilde{\textbf{E}},\widetilde{\textbf{K}} \in \mathbb{R}^{r \times r}$, $\widetilde{\textbf{B}} \in \mathbb{R}^{r \times m}$, with $r \ll n$.
The properties of the original system, i.e., 
\begin{equation} \label{eq:full_constraints}
\begin{aligned}
\textbf{M} &\succ 0,  \quad \textbf{K} \succ 0, \\
\textbf{D} &\succeq 0, \quad \textbf{G} = - \textbf{G}^{\mathsf{T}}
\end{aligned}
\end{equation}
 should be preserved in the reduced-model to guarantee interpretable results for different load cases.

\section{Reduced-order inference with embedded parametrization}\label{sec:popinf}

As previously stated, in this work we pursue two essential objectives:
\begin{itemize}
	\item  Obtaining the full physical system representation in a low-dimension  $\widetilde{\boldsymbol{\Sigma}} = (\widetilde{\textbf{M}}, \widetilde{\textbf{E}}, \widetilde{\textbf{K}}, \widetilde{\textbf{\textbf{B}}})$, where $\widetilde{\textbf{E}} = \widetilde{\textbf{D}} + \Omega \widetilde{\textbf{G}}$, from the available data, and
	\item  Preservation of the SPD properties for the inferred system operators $\widetilde{\textbf{M}} \succ 0$, $\widetilde{\textbf{K}} \succ 0$, $\widetilde{\textbf{D}} \succeq 0$, and skew-symmetric properties for the inferred gyroscopic term $\widetilde{\textbf{G}} =  - \smash{\widetilde{\textbf{G}}}^{\mathsf{T}} $, if the system is subjected to rotation.
\end{itemize}

The seed idea is that the reduced-system satisfies the equations of motion for low-dimensional data approximation $\textbf{x}_r$, $\dot{\textbf{x}}_r$, and $\ddot{\textbf{x}}_r \in \mathbb{R}^{r}$, $r \ll n$:

\begin{equation}
\label{eq:red_noG}
\widetilde{\textbf{M}}\ddot{\mathbf{x}}_r(t) + \widetilde{\textbf{E}}\dot{\mathbf{x}}_r(t) + \widetilde{\textbf{K}}\mathbf{x}_r(t) = \widetilde{\textbf{B}} \mathbf{u}(t).
\end{equation}
We aim to develop a data-driven framework based on the ideas, presented in \cite{morPehW16, filanova2022operator}, to construct a low-dimensional second-order dynamical system, capturing the dynamics present in the given data. 

Before proceeding with the method description, we clarify the assumptions concerning the data for reduced-order modeling in \Cref{subsec:problem}, and introduce the necessary notation. The following \Cref{subsec:method} describes the reduced-order inference methodology, where the parametrization of the unknown variables is used in the nonlinear optimization problem to preserve the SPD and skew-symmetric structure of the system matrices.

\subsection{Available data}\label{subsec:problem}

An important part of choosing the appropriate nonintrusive reduction method is the analysis of the accessible data. In this paper, we consider the case where the input signal $\textbf{u}(t)$ and the corresponding system response $\textbf{x}(t)$ of \eqref{eq:original_noG} are known from a simulation over a time interval $t \in [0, \; T ]$. The derivatives of the state vector are often approximated by some numerical scheme \cite{Kramer2024}. Often, the approximation of the derivative snapshots causes additional inaccuracies and drastically affects the results. One way to avoid the calculation of the derivatives is to embed the integration scheme into the optimization process \cite{peh_rollout}.
On the other hand, the software  for the simulation of mechanical systems usually uses the Newmark family of time integration algorithms, which provide the derivative information for each time point, satisfying the governing system of equations \cite{GarciadeJalon1994}. Since the accelerations and velocities can be provided by second-order time-integrators, as shown in the algorithm \cite{HilHT77}, we will use that given derivative information for the further optimization. 

Moreover, we assume that the input operator $\textbf{B}$ is unknown. One of the possible reasons for the latter assumption is the difficulty to get the correct $\textbf{B}$ operator if there exist a software-specific difference in the node numbering for the solution process and for displaying the results. Finally, we consider the rotation speed $\Omega$ to be known, as it is a part of an external input to the system.

Now, we introduce the notation used for reduced-order modeling. First, we denote the snapshot state and input matrices as $\textbf{X}$ and $\textbf{U}$, respectively. These matrices contain the inputs and the state snapshots at pre-defined time steps into the corresponding columns:
\begin{equation}
\label{mat:snap}
\textbf{U} =
\begin{bmatrix}
| &  \dots & | \\
\mathbf{u}(t_1) & \dots & \mathbf{u}(t_N) \\
| & \dots & |
\end{bmatrix} \in \mathbb{R}^{m \times N}, \quad
\textbf{X} = \begin{bmatrix}
| &  \dots & | \\
\mathbf{x}(t_1) & \dots & \mathbf{x}(t_N) \\
| & \dots & |
\end{bmatrix} \in \mathbb{R}^{n \times N}.
\end{equation} The respective derivative snapshot matrices are assembled analogously:
\begin{equation} \label{mat:derivative}
\begin{aligned}
\displaystyle \dot{\textbf{X}} &= \begin{bmatrix}
| &  \dots & | \\
\dot{\mathbf{x}}(t_1) & \dots & \dot{\mathbf{x}}(t_N) \\
| & \dots & |
\end{bmatrix} \in \mathbb{R}^{n \times N}, \quad
\displaystyle \ddot{\textbf{X}} &= \begin{bmatrix}
| &  \dots & | \\
\ddot{\mathbf{x}}(t_1) & \dots & \ddot{\mathbf{x}}(t_N) \\
| & \dots & |
\end{bmatrix} \in \mathbb{R}^{n \times N}.
\end{aligned}
\end{equation}

The final necessary procedure for data preparation is to compress the data \eqref{mat:snap} and \eqref{mat:derivative}, such that we can obtain a low-dimensional surrogate model that satisfies \eqref{eq:red_noG}. The reduced snapshots can be constructed by projecting the original solution snapshots onto the proper orthogonal decomposition (POD) subspace:

\begin{equation}\label{eq:low_dimensionaldata}
\textbf{X}_r = \textbf{V}_r^{\mathsf{T}} \textbf{X}.
\end{equation}
The projection basis $\textbf{V}_r \in \mathbb{R}^{n \times r}$ is the truncated left-singular vector matrix resulting from the singular value decomposition of the snapshot matrix:

\begin{equation}
\label{eq:svd}
\textbf{X} = \textbf{V} \boldsymbol{\Sigma} \textbf{W}^{\mathsf{T}}.
\end{equation} 

The next subsection contains the details of the operator inference methodology with the embedded parametrization.

\subsection{Methodology description} \label{subsec:method}

For the sake of brevity, in this subsection we will consider the general form of the second-order system given in \eqref{eq:original}, including the gyroscopic term, which just vanishes from the optimization process, if there are no rotation degrees of freedom in the system.  

%Following the operator inference approach in its second-order reformulation, presented in \cite{filanova2022operator}, we state an optimization problem:
%
%\begin{align}\label{eq:fopinf}
%\min_{\substack{\widetilde{\textbf{M}}, \, \widetilde{\textbf{D}}, \\ \widetilde{\textbf{G}}, \, \widetilde{\textbf{K}}, \\ \widetilde{\textbf{B}}}} \sum_{i = 0}^{N} \| \ddot{\mathbf{x}}_r(t_i) + \widetilde{\textbf{M}}^{-1} \widetilde{\textbf{D}} \dot{\mathbf{x}}_r(t_i)  + \Omega \widetilde{\textbf{M}}^{-1} \widetilde{\textbf{G}} \dot{\mathbf{x}}_r(t_i) + \widetilde{\textbf{M}}^{-1} \widetilde{\textbf{K}} \mathbf{x}(t) - \widetilde{\textbf{M}}^{-1} \widetilde{\textbf{B}}  \mathbf{u}(t_i) \| _F^2.
%\end{align}

%Now, instead of considering the terms $\widetilde{\textbf{M}}^{-1} \widetilde{\textbf{K}}$, etc. as one coupled unknown matrix-parameter and solving \eqref{eq:fopinf} as a linear least-squares problem, we propose to optimize a function over each unknown matrix-parameter $\widetilde{\textbf{M}},\widetilde{\textbf{D}},\widetilde{\textbf{G}}, \widetilde{\textbf{K}}, \; \text{and} \; \widetilde{\textbf{B}}$.

%The minimization problem \eqref{eq:fopinf} can be reformulated to 
First, we formulate a loss function $\mathcal{F}$, which compares the true compressed displacements $\bX_r$ \eqref{eq:low_dimensionaldata}   and the predicted displacements $\bX_p$:

\begin{equation} \label{eq:mse_func}
\mathcal{F}  = \frac{1}{N} \; \| \bX_r - \bX_p \|_F^2.
\end{equation} The predicted displacements $\bX_p$ are calculated as

\begin{equation}\label{eq:xp}
\bX_p = \smash{\widetilde{\textbf{K}}}^{-1} \widetilde{\textbf{B}}  \mathbf{U} - \smash{\widetilde{\textbf{K}}}^{-1} \tilde{\textbf{M}} \ddot{\mathbf{X}}_r - \smash{\widetilde{\textbf{K}}}^{-1} \left( \widetilde{\textbf{D}} + \Omega \widetilde{\textbf{G}} \right) \dot{\mathbf{X}}_r. 
\end{equation} Note that in the optimization process, the inverse of the stiffness matrix is an independent variable, i.e., $\smash{\widetilde{\textbf{K}}}^{-1} \equiv  \widetilde{\textbf{K}}_{\text{inv}}$ and there is no matrix inversion during the calculation of $\bX_p$.

Thus, in order to obtain the full system representation, we aim to minimize the function $\mathcal{F} (\widetilde{\textbf{M}}, \widetilde{\textbf{K}}_{\text{inv}},\widetilde{\textbf{D}},\widetilde{\textbf{G}}, \widetilde{\textbf{B}})$ defined in \cref{eq:mse_func} and \cref{eq:xp}. The resulting nonlinear optimization problem allows us to find the full physical system representation in reduced dimension.

Now, it remains us to ensure that the solution of the optimization problem described above preserves the SPD and skew-symmetric properties of the corresponding operators \eqref{eq:full_constraints}. For this purpose, we consider the following parametrization of the unknown system matrices:

\begin{equation}
\begin{aligned}
\widetilde{\textbf{K}}_{\text{inv}}  & = \widetilde{\textbf{K}}_c \smash{\widetilde{\textbf{K}}}_c^{\mathsf{T}} , \quad 
\widetilde{\textbf{M}}  = \widetilde{\textbf{M}}_c \smash{\widetilde{\textbf{M}}}_c^{\mathsf{T}},  \\
\widetilde{\textbf{D}} & = \widetilde{\textbf{D}}_c \smash{\widetilde{\textbf{D}}}_c^{\mathsf{T}}, \quad
\widetilde{\textbf{G}} = \widetilde{\textbf{G}}_c - \smash{\widetilde{\textbf{G}}}_c^{\mathsf{T}}. \label{eq:decomposed_vars2}
\end{aligned}
\end{equation} Including the parametrization  \eqref{eq:decomposed_vars2} into \eqref{eq:xp} allows us to avoid using the linear matrix inequality constraints to ensure the SPD and skew-symmetric system matrices. Thus, the final optimization problem appears as follows:

\begin{align}\label{eq:popinf_opt}
	\min_{\substack{\widetilde{\textbf{M}}_c, \, \widetilde{\textbf{D}}_c, \\ \widetilde{\textbf{G}}_c, \, \widetilde{\textbf{K}}_c, \\ \widetilde{\textbf{B}}}} \| \mathbf{X}_r - \widetilde{\textbf{K}}_c \smash{\widetilde{\textbf{K}}}_c^{\mathsf{T}}  \widetilde{\textbf{B}}  \mathbf{U} - \widetilde{\textbf{K}}_c \smash{\widetilde{\textbf{K}}}_c^{\mathsf{T}}  \widetilde{\textbf{M}}_c \smash{\widetilde{\textbf{M}}}_c^{\mathsf{T}} \ddot{\mathbf{X}}_r - \widetilde{\textbf{K}}_c \smash{\widetilde{\textbf{K}}}_c^{\mathsf{T}}  \left( \widetilde{\textbf{D}}_c \smash{\widetilde{\textbf{D}}}_c^{\mathsf{T}}  + \Omega (\widetilde{\textbf{G}}_c - \smash{\widetilde{\textbf{G}}}_c^{\mathsf{T}}) \right) \dot{\mathbf{X}}_r \|_F^2.
\end{align} The optimization problem \eqref{eq:popinf_opt} is nonlinear and can be solved using the efficient algorithms based on automatic differentiation, e.g., implemented in PyTorch framework \cite{pytorch}. 

In the following subsection, we will discuss some of the implementation details that we used for the numerical examples and that might strongly influence the solution.

\subsection{Implementation details}
\label{subsec:implementation}
The practical solution of the inference problem \eqref{eq:popinf_opt} is implemented in the PyTorch framework using the Adam optimizer \cite{adam}. In order to reach a reasonable (sub-)optimum efficiently, many hyperparameters have to be tuned. 

%In the following, we will discuss some optimizer-specific, as well as simulation-specific parameters, chosen for the numerical examples in \Cref{sec:numres}. 

The choice of the simulation-specific parameters, such as the time-stepping scheme, was based depending on the respective model and is described in the \Cref{sec:numres}. 

For each numerical example, the number of epochs was set to 36 000. The learning rate for each iteration while minimizing the loss function was specified according to the cyclic learning rate policy, as it showed the best performance. Unlike the constant or monotonically decreasing learning rate, it allows the variation between specified boundary values. For the numerical examples in this paper, we chose the lower bound for the learning rate to be $5 \cdot 10^{-6}$, and the upper bound $1 \cdot 10^{-3}$.

Additionally, we would like to emphasize that the system operators typically have very different orders of magnitude for structural mechanics simulations. For metal parts, the stiffness matrix can be of order $10^9$ or more, while the mass matrix is of order $10^{-2} - 10^1$, depending on the geometric size, etc. This drastic difference causes problems in solving the optimization problem, and affects the accuracy. To avoid these issues, we normalize the snapshot matrices by their corresponding norm, denoted as $\alpha_X, \alpha_V, \alpha_A$, and $\alpha_U$:
\begin{align}
	\textbf{X}_r \vcentcolon =  \frac{1}{\alpha_X} \textbf{X}_r,  \quad \dot{\textbf{X}}_r := \frac{1}{\alpha_V} \dot{\textbf{X}}_r, \quad
	\ddot{\textbf{X}}_r : = \frac{1}{\alpha_A} \ddot{\textbf{X}}_r, \quad \textbf{U} := \frac{1}{\alpha_U} \textbf{U}.
\end{align}
 The corresponding operators are then scaled to the correct magnitude after the optimum is found:
 \begin{align}
 \widetilde{\textbf{K}} := \frac{1}{\alpha_X} \widetilde{\textbf{K}},  \quad \widetilde{\textbf{D}} := \frac{1}{\alpha_V} \widetilde{\textbf{D}},
 \quad \widetilde{\textbf{G}} := \frac{1}{\alpha_V} \widetilde{\textbf{G}}, 
\quad  \widetilde{\textbf{M}} := \frac{1}{\alpha_A} \widetilde{\textbf{M}}, \quad \widetilde{\textbf{B}} := \frac{1}{\alpha_U} \widetilde{\textbf{B}}.
 \end{align}

To solve the optimization problem \eqref{eq:popinf_opt}, initial values for the unknown coefficients had to be specified. The initial guesses can have a significant impact on the results. In this paper, we use homogeneous initial values for the gyroscopic matrix, ones for the control matrix, and random values with uniform distribution on the interval $[0,1)$ for the other system matrices.

%For every operation the respective derivative is calculated, such that the gradient calculation can be done as efficient as the calculation of the function itself.

% Usage & trainings aspects: influence of r, data scaling, learning rate and so on?

\section{Numerical results}\label{sec:numres}
In the following section, we show the numerical results of the parametrized inference method (\texttt{p-OpInf}) for three examples. Additionally, we provide the comparison of the ROMs identified by \texttt{p-OpInf} and the constrained force-informed operator inference method (\texttt{fi-OpInf}) described in \cite{filanova2022operator}. The latter approach is implemented in Python using the CVXPY package~\cite{diamond2016cvxpy}. 

The numerical experiments are designed as follows. First, we run simulations and collect the snapshots to construct the reduced-order models. Since the snapshot data should contain enough "dynamics", our goal is to stimulate a wide range of frequencies and obtain a "rich" system response. Therefore, we use the chirp input signal $\mathrm{u}_{\text{m}}(t)$ with amplitude $\mathrm{A_m}$ for the reduced-order modeling: 
\begin{equation} \label{eq:chirp_signal}
\begin{split}
 \mathrm{u}_{\text{m}}(t) & = \mathrm{A_m} \cos (\phi(t)), \\
\text{with} \quad \phi(t) &  = \phi_0 + 2 \pi \int_0^{t} f(\tau) d \tau.
\end{split}
\end{equation} In \eqref{eq:chirp_signal}, $\phi_0$ is the initial phase shift, and $f(t)$ is the input frequency that changes over time as follows:
\begin{equation} \label{eq:chirp}
	f(t) = \frac{f_1 - f_0}{\Delta t}t + f_0,
\end{equation} where $\Delta t$ is the time interval to sweep from the initial frequency $f_0$ to the target frequency $f_1$. For the following numerical examples, we choose the frequency range $[f_0 , f_1]$ depending on model specifics such as eigenfrequencies, but in general we do not provide guidelines for the choice of loading conditions to generate snapshot data.

The numerical integration of the high-fidelity and reduced-order models for all the examples is performed using the Python implementation of the Newmark-beta algorithm \cite{New59}. Since this method is widely used in commercial finite element tools and provides the velocity and acceleration data for each time step, for the reduced-order modeling we will use the derivative snapshots directly from the integrator. 

The quality of the results is measured by means of relative error $\epsilon_{\mathrm{err}}$. For each time step $t_i$ we divide the absolute error by the maximum norm of the original trajectory vector to avoid the distorted error behavior when the displacements are close to zero:

\begin{equation} \label{eq:rel_error}
\epsilon_{\mathrm{err}} (t_i) = \frac{\left\| \mathbf{x}(t_i) - \hat{\mathbf{x}}(t_i)\right\|_2}{\max_{t\in [t_1, \; t_N]} \left\| \mathbf{x}(t) \right\|_2},
\end{equation} where $\hat{\mathbf{x}} (t)$ is the trajectory vector of the obtained reduced-order model. We validate the results by simulating the ROMs using the same loading conditions as for the snapshots generation, i.e., input signal, force amplitude, time step, etc. Different test loading conditions are proposed to further analyze the ROM behavior. 

%All simulations were performed on Desktop PC Dell Precision 3660 Tower Intel \textsuperscript{\textregistered} Core \textsuperscript{\texttrademark} i5-12600K, 31GB RAM. T

\subsection{Structural beam}
\label{subsec:beam}
\begin{figure}[h!]
	\centering
	\begin{subfigure}[b]{0.47\textwidth}
\includegraphics[scale = 0.25]{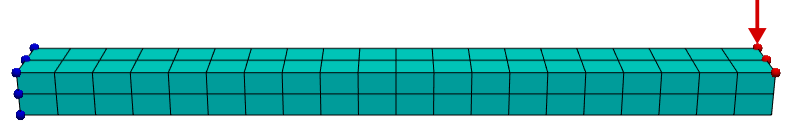}
\caption{Cantilever beam.}
\label{fig:ex1}
	\end{subfigure}
	\hfill
	\begin{subfigure}[b]{0.45\textwidth}
	\includegraphics[scale = 0.25]{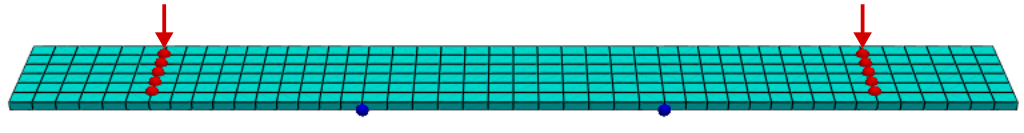}
\caption{Overhanging beam.}
\label{fig:ex2}
	\end{subfigure}
\caption{Finite element meshes of the beam models. The nodes with applied boundary conditions are marked with blue color, the positions of the external force are marked with red color.}
\label{fig:beams}
\end{figure}

First, we consider two steal beam models, which are common examples of structural mechanics problems. \Cref{fig:ex1} shows a cantilever beam model of dimension 540, while a slightly larger model with 6225 degrees of freedom (DOFs) is shown in \Cref{fig:ex2}. The nodes with applied boundary and load conditions are marked in \Cref{fig:beams} with blue and red color, respectively. The nodes with boundary conditions are fixed in all directions.
%The finite element meshes and system matrices are imported from Ansys \textsuperscript{\textregistered} Academic Student Mechanical, Release 24.2. 

  The cantilever beam is a common structural element in engineering practice, typically loaded at the unsupported end. In the following, we assume the force to act on three nodes, as shown in \Cref{fig:ex1}, having the load direction perpendicular to the beam axis. The overhanging beam can typically be used for the four-point bending tests, correspondingly the force is located at two unsupported ends and applied as shown in \Cref{fig:ex2}.

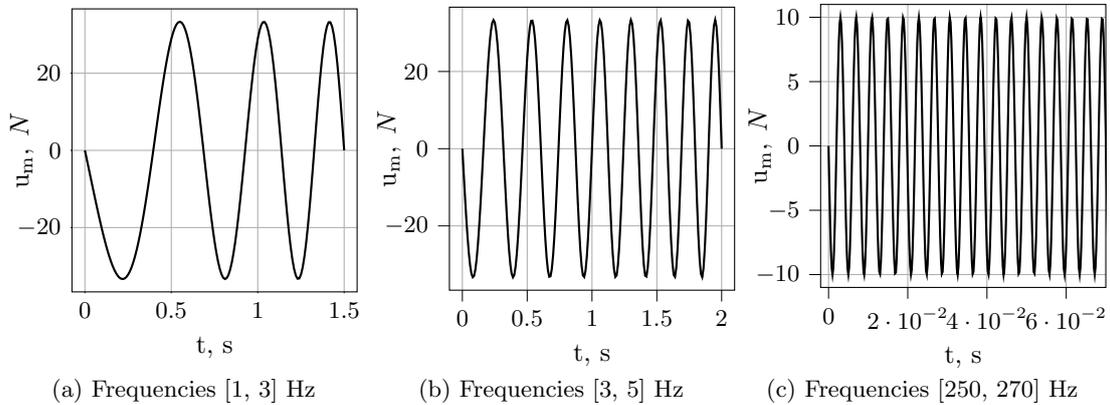
\begin{figure}[h!]
	\centering
	\begin{subfigure}[b]{0.32\textwidth}
		% This file was created with tikzplotlib v0.10.1.
\begin{tikzpicture}

\definecolor{darkgray176}{RGB}{176,176,176}

\begin{axis}[
width=\plotwtripple,
height=\plothtripple,
tick align=outside,
tick pos=left,
major tick length = 0.5pt,
x grid style={darkgray176},
xlabel={t, s},
xmajorgrids,
xmin=-0.075, xmax=1.575,
xtick style={color=black},
y grid style={darkgray176},
ylabel={\(\displaystyle \mathrm{u_m}, \; N\)},
ylabel style={yshift=-13pt},
ymajorgrids,
ymin=-36.6630257938804, ymax=36.6645785712027,
ytick style={color=black},
xticklabel style={font = \small},
yticklabel style={font = \small}
]
\addplot [thick, black]
table {%
0 2.04107799857892e-15
0.01 -2.10695221556343
0.02 -4.23317869987788
0.03 -6.36943699850649
0.04 -8.50585866521812
0.05 -10.6319769766023
0.06 -12.7367622078059
0.07 -14.8086648441126
0.08 -16.8356670537784
0.09 -18.8053426885529
0.1 -20.704926009277
0.11 -22.5213892545786
0.12 -24.2415290808132
0.13 -25.8520618010039
0.14 -27.3397272397418
0.15 -28.6914009001315
0.16 -29.8942140083929
0.17 -30.9356808623831
0.18 -31.8038327630156
0.19 -32.4873576535145
0.2 -32.9757444320996
0.21 -33.2594307407554
0.22 -33.3299528681948
0.23 -33.1800962412425
0.24 -32.8040448181848
0.25 -32.1975275429689
0.26 -31.3579598735616
0.27 -30.2845782645938
0.28 -28.9785653671448
0.29 -27.4431636108418
0.3 -25.683774759193
0.31 -23.7080429821303
0.32 -21.5259189740343
0.33 -19.1497026649086
0.34 -16.5940621306145
0.35 -13.87602640868
0.36 -11.0149500723826
0.37 -8.03244761037129
0.38 -4.95229590433118
0.39 -1.80030339475495
0.4 1.39585512430665
0.41 4.60683873653354
0.42 7.80187585917126
0.43 10.9490346149937
0.44 14.0155272577832
0.45 16.968047191679
0.46 19.7731364906109
0.47 22.3975811599443
0.48 24.8088306969856
0.49 26.9754378101189
0.5 28.8675134594813
0.51 30.4571916979696
0.52 31.7190981340215
0.53 32.6308152221496
0.54 33.1733370296978
0.55 33.3315056455171
0.56 33.0944210054701
0.57 32.4558156282421
0.58 31.4143856000225
0.59 29.9740691347633
0.6 28.1442641834005
0.61 25.9399768845926
0.62 23.3818931531075
0.63 20.4963663993654
0.64 17.3153152711296
0.65 13.87602640868
0.66 10.2208585066627
0.67 6.39684547331066
0.68 2.45519816005235
0.69 -1.54929401446377
0.7 -5.55895822387013
0.71 -9.5140247373879
0.72 -13.3534971230641
0.73 -17.0160962551537
0.74 -20.4412659743562
0.75 -23.5702260395516
0.76 -26.3470558984769
0.77 -28.7197908493254
0.78 -30.6415104249972
0.79 -32.071397367454
0.8 -32.9757444320996
0.81 -33.3288855305348
0.82 -33.1140274403865
0.83 -32.3239585341165
0.84 -30.9616117487733
0.85 -29.040460370673
0.86 -26.5847271673253
0.87 -23.6293899751844
0.88 -20.2199700432541
0.89 -16.4120932204427
0.9 -12.2708184228226
0.91 -7.8697326713021
0.92 -3.28981727753812
0.93 1.38190461584284
0.94 6.05392307607017
0.95 10.6319769766024
0.96 15.0209422043424
0.97 19.1268388129231
0.98 22.858918104427
0.99 26.1317858112908
1 28.8675134594813
1.01 30.9976867988888
1.02 32.465338055309
1.03 33.22670784517
1.04 33.2527830348595
1.05 32.5305587312916
1.06 31.0639760380491
1.07 28.8744922437435
1.08 26.0012467255408
1.09 22.5007940037016
1.1 18.4463849747781
1.11 13.9267882315632
1.12 9.04465534334144
1.13 3.91444676426422
1.14 -1.34005164543181
1.15 -6.58857801263764
1.16 -11.6975530175555
1.17 -16.5334782624588
1.18 -20.9664935236696
1.19 -24.8740019457801
1.2 -28.1442641834005
1.21 -30.6798568282504
1.22 -32.4008875581757
1.23 -33.2478596345668
1.24 -33.1840819088555
1.25 -32.197527542969
1.26 -30.3020552664848
1.27 -27.5379211340827
1.28 -23.9715262360505
1.29 -19.6943663472917
1.3 -14.8211726394977
1.31 -9.48725776029791
1.32 -3.84510811361026
1.33 1.93970974063269
1.34 7.69323390012532
1.35 13.2382630211595
1.36 18.3998400059899
1.37 23.010900283322
1.38 26.9179058456747
1.39 29.986275709351
1.4 32.1054188932553
1.41 33.1931780029156
1.42 33.199500495661
1.43 32.1091709055779
1.44 29.9434606803902
1.45 26.7605825063705
1.46 22.6548724599345
1.47 17.7546651195566
1.48 12.2188727077008
1.49 6.23232795458053
1.5 -3.2677880665149e-14
};
\end{axis}

\end{tikzpicture}
		\vspace{-0.5cm}
		\caption{Frequencies [1, 3] \si{\hertz}}
	\end{subfigure}
%\hfill
\hspace{-0.2cm}
	\begin{subfigure}[b]{0.32\textwidth}
		% This file was created with tikzplotlib v0.10.1.
\begin{tikzpicture}

\definecolor{darkgray176}{RGB}{176,176,176}

\begin{axis}[
width=\plotwtripple,
height=\plothtripple,
tick align=outside,
tick pos=left,
x grid style={darkgray176},
xlabel={t, s},
xmajorgrids,
xmin=-0.1, xmax=2.1,
xtick style={color=black},
y grid style={darkgray176},
ylabel={\(\displaystyle \mathrm{u_m}, \; N\)},
ylabel style={yshift=-13pt},
ymajorgrids,
ymin=-36.6558405108969, ymax=36.6651585399278,
ytick style={color=black},
xticklabel style={font = \small},
yticklabel style={font = \small}
]
\addplot [thick, black]
table {%
0 2.04107799857892e-15
0.01 -6.25632999916038
0.02 -12.309755110291
0.03 -17.940397698847
0.04 -22.9400879514852
0.05 -27.1202816500262
0.06 -30.3194844172945
0.07 -32.4098844463447
0.08 -33.3029132608633
0.09 -32.9534833729989
0.1 -31.3626922984742
0.11 -28.5788332644354
0.12 -24.6966126317586
0.13 -19.8545406294402
0.14 -14.2305331430775
0.15 -8.03583536110129
0.16 -1.5074501694311
0.17 5.10067772896207
0.18 11.527409950043
0.19 17.5147317280598
0.2 22.8182368642896
0.21 27.2173553415224
0.22 30.5248919000567
0.23 32.5954550809029
0.24 33.3323858557994
0.25 32.6928426801077
0.26 30.6907645101733
0.27 27.397513092824
0.28 22.9400879514851
0.29 17.4969085416047
0.3 11.2912640081764
0.31 4.58263737520717
0.32 -2.3437869235101
0.33 -9.18903033374002
0.34 -15.6531613650047
0.35 -21.4485288244593
0.36 -26.3128062326085
0.37 -30.0212523513902
0.38 -32.3976053023848
0.39 -33.3230678267685
0.4 -32.742908357623
0.41 -30.6702949430872
0.42 -27.1870935695953
0.43 -22.4414947418517
0.44 -16.6424768293183
0.45 -10.0512653314832
0.46 -2.9700968422721
0.47 4.27126223592711
0.48 11.3306666040688
0.49 17.8697340657813
0.5 23.5702260395516
0.51 28.1498739595201
0.52 31.3768565525769
0.53 33.0821702181148
0.54 33.1692116538936
0.55 31.6200064875168
0.56 28.4976659584984
0.57 23.944829659147
0.58 18.1780482273512
0.59 11.4782664744278
0.6 4.17777445214353
0.61 -3.35580943246611
0.62 -10.7377838305989
0.63 -17.5859552581712
0.64 -23.5405882174572
0.65 -28.2836738327168
0.66 -31.5564749222476
0.67 -33.1743601887097
0.68 -33.038050658613
0.69 -31.1405642085729
0.7 -27.5693524758187
0.71 -22.5033692739785
0.72 -16.2050789110616
0.73 -9.00769263141776
0.74 -1.29819655889632
0.75 6.50301073387084
0.76 13.9648339488345
0.77 20.6693439531166
0.78 26.2354745464155
0.79 30.3412035233011
0.8 32.742908357623
0.81 33.2907040867499
0.82 31.9387622174943
0.83 28.7498668257055
0.84 23.8937752683578
0.85 17.6393000324731
0.86 10.3403960978765
0.87 2.41690359403947
0.88 -5.66906362029641
0.89 -13.4397982440237
0.9 -20.4302351217659
0.91 -26.216083321793
0.92 -30.4401463402829
0.93 -32.8351754807543
0.94 -33.2417590455414
0.95 -31.6200064875168
0.96 -28.0541303718739
0.97 -22.7494420725296
0.98 -16.0217364854271
0.99 -8.27951951743898
1 -8.98280646253588e-14
1.01 8.29980547539455
1.02 16.0951497628604
1.03 22.8868492385247
1.04 28.2336868981751
1.05 31.7818057439667
1.06 33.2891009726431
1.07 32.6428823268581
1.08 29.8694647523332
1.09 25.134835700667
1.1 18.7361125950711
1.11 11.0841092963293
1.12 2.677936671113
1.13 -5.9268715246188
1.14 -14.1547305491403
1.15 -21.4485288244592
1.16 -27.307737422748
1.17 -31.3234852989543
1.18 -33.2081073222709
1.19 -32.8169766942727
1.2 -30.1609017488674
1.21 -25.4079731588076
1.22 -18.8744539972698
1.23 -11.0050658269173
1.24 -2.34378692351027
1.25 6.50301073387082
1.26 14.9086542698624
1.27 22.270610045009
1.28 28.0541303718738
1.29 31.8319271996493
1.3 33.3168853455244
1.31 32.3852462205543
1.32 29.0883156301529
1.33 23.6515379910086
1.34 16.4606819125001
1.35 8.0358353611011
1.36 -1.00515725402856
1.37 -9.99134021813624
1.38 -18.2482127938918
1.39 -25.1485872003524
1.4 -30.1609017488673
1.41 -32.8912329829455
1.42 -33.1156235505847
1.43 -30.7999903487886
1.44 -26.1057603608681
1.45 -19.3804394773278
1.46 -11.1334759427783
1.47 -1.99894727364588
1.48 7.3123147335574
1.49 16.0676316728809
1.5 23.5702260395516
1.51 29.2152661392571
1.52 32.5396706065184
1.53 33.2615067784618
1.54 31.3055402471889
1.55 26.8125211756695
1.56 20.1310534477092
1.57 11.7922901580217
1.58 2.46912265172859
1.59 -7.07713724302129
1.6 -16.0584558033906
1.61 -23.7252128697031
1.62 -29.4294011218753
1.63 -32.6805267211588
1.64 -33.1893146167381
1.65 -30.8952198546706
1.66 -25.9750155646015
1.67 -18.831273372875
1.68 -10.0612493840532
1.69 -0.408396826924782
1.7 9.29970353464087
1.71 18.2219145653411
1.72 25.576649329506
1.73 30.7111583508488
1.74 33.1608037727443
1.75 32.6928426801077
1.76 29.3304691443203
1.77 23.354511470844
1.78 15.2821144289457
1.79 5.82379134450951
1.8 -4.17777445214369
1.81 -13.8220492921036
1.82 -22.2316239538237
1.83 -28.6325923358387
1.84 -32.4269408946171
1.85 -33.2500932138757
1.86 -31.0079450200213
1.87 -25.8894812253265
1.88 -18.3532433823464
1.89 -9.08832307153644
1.9 1.04702530260444
1.91 11.1038592436042
1.92 20.1310534477092
1.93 27.2656344608788
1.94 31.8163505817987
1.95 33.3323052548263
1.96 31.6497073081463
1.97 26.9117279716512
1.98 19.5589382575003
1.99 10.2906065304671
2 8.14367930642941e-15
};
\end{axis}

\end{tikzpicture}
		\vspace{-0.61cm}
		\caption{Frequencies [3, 5] \si{\hertz}}
	\end{subfigure}
%\hfill
\hspace{-0.15cm}
	\begin{subfigure}[b]{0.32\textwidth}
	% This file was created with tikzplotlib v0.10.1.
\begin{tikzpicture}

\definecolor{darkgray176}{RGB}{176,176,176}

\begin{axis}[
width=\plotwtripple,
height=\plothtripple,
tick align=outside,
tick pos=left,
x grid style={darkgray176},
xlabel={t, s},
xmin=-0.002, xmax=0.1,
xtick style={color=black},
y grid style={darkgray176},
ymajorgrids,
xmajorgrids,
ylabel={\(\displaystyle \mathrm{u_m}, \; N\)},
ylabel style={yshift=-13pt},
ymin=-10.999894902704, ymax=10.9999820708876,
ytick style={color=black},
scaled x ticks=false,
xticklabel style={font = \small},
yticklabel style={font = \small},
xtick distance={0.05}
%xticklabels={
%	\(\displaystyle {0.0}\),
%	\(\displaystyle {0.05}\),
%	\(\displaystyle {0.1}\),
%}
]
\addplot [thick, black]
table[x = time, y = um,  col sep=comma] {pics/beam4/csv/ex2_um_250_270.csv};
\end{axis}

\end{tikzpicture}
	\vspace*{-0.64cm}
	\caption{Frequencies [250, 270] \si{\hertz}}
\end{subfigure}
	\caption{Input signals $\mathrm{u_m}(t)$ according to \eqref{eq:chirp} with different frequency ranges $[f_0, f_1]$\si{\hertz}.}
	\label{fig:beam_input}
\end{figure}

The governing system of equations for both beam models corresponds to \eqref{eq:original_noG}. For the data collection, we apply the harmonic input signal $\mathrm{u_m} (t)$ with variable frequency from \eqref{eq:chirp}, as shown in \Cref{fig:beam_input}.

 The frequencies $f_0$ and $f_1$ can be chosen arbitrarily. However, we decided to focus on the system behavior around the first three eigenfrequencies for the cantilever beam, which are [0.98, 1.63, 6.11] \si{\hertz} , and the first two eigenfrequencies for the overhanging beam, which are [268.28 , 278.16] \si{\hertz}. We integrate the system for three different frequency ranges and form the corresponding snapshot matrices $\bX$, as explained in \Cref{subsec:problem}. The singular value decay of the different snapshot matrices is plotted in \Cref{fig:beam_svds} for both beam models.
%[268.28 , 278.16, 1674.79] \si{\hertz}
\begin{figure}[h!]
	\centering
	\begin{subfigure}[b]{0.49\textwidth}
		% This file was created with tikzplotlib v0.10.1.
\begin{tikzpicture}

\definecolor{darkgray176}{RGB}{176,176,176}
\definecolor{darkred}{RGB}{139,0,0}
\definecolor{dimgray}{RGB}{105,105,105}
\definecolor{lightgray204}{RGB}{204,204,204}
\definecolor{slategray}{RGB}{112,128,144}

\begin{axis}[
width=\plotw,
height=\ploth,
legend cell align={left},
legend style={fill opacity=0.8, draw opacity=1, text opacity=1, draw=lightgray204},
log basis y={10},
tick align=outside,
tick pos=left,
x grid style={darkgray176},
xlabel={r},
xmajorgrids,
xmin=-3, xmax=81,
xtick style={color=black},
y grid style={darkgray176},
ylabel={\(\displaystyle \sigma / \sigma_{max}\)},
ymajorgrids,
ymin=1.41253754462275e-18, ymax=7.07945784384138,
ymode=log,
ytick style={color=black},
ytick={1e-21,1e-18,1e-15,1e-12,1e-09,1e-06,0.001,1,1000,1000000},
yticklabels={
  \(\displaystyle {10^{-21}}\),
  \(\displaystyle {10^{-18}}\),
  \(\displaystyle {10^{-15}}\),
  \(\displaystyle {10^{-12}}\),
  \(\displaystyle {10^{-9}}\),
  \(\displaystyle {10^{-6}}\),
  \(\displaystyle {10^{-3}}\),
  \(\displaystyle {10^{0}}\),
  \(\displaystyle {10^{3}}\),
  \(\displaystyle {10^{6}}\)
},
xticklabel style={font = \small},
yticklabel style={font = \small}
]
\addplot [semithick, black, mark=o, mark size=2, mark options={solid}, only marks]
table {%
0 1
1 0.015312047581746
2 0.000310147898712378
3 2.75959407005178e-05
4 4.82097124296808e-06
5 1.22441070595319e-06
6 4.31622954063824e-07
7 1.78100661919927e-07
8 8.67743986500293e-08
9 4.66551327401354e-08
10 2.89110017330974e-08
11 1.59221721255257e-08
12 1.01058934252917e-08
13 5.80737521174464e-09
14 5.69132986391671e-09
15 1.43654960925018e-09
16 7.3507025880727e-10
17 6.50595939404898e-10
18 4.45793246771158e-10
19 7.20553994530257e-12
20 1.82459775205363e-13
21 1.39011399807396e-13
22 5.28320010923526e-14
23 7.82954201858568e-15
24 5.28560369428251e-15
25 3.30084912650579e-15
26 7.04535720039237e-16
27 5.22727513663356e-16
28 2.30687638210182e-16
29 1.74731601905915e-16
30 1.39049680030563e-16
31 1.06768806806476e-16
32 9.59500405313704e-17
33 9.59500405313704e-17
34 9.59500405313704e-17
35 9.59500405313704e-17
36 9.59500405313704e-17
37 9.59500405313704e-17
38 9.59500405313704e-17
39 9.59500405313704e-17
40 9.59500405313704e-17
41 9.59500405313704e-17
42 9.59500405313704e-17
43 9.59500405313704e-17
44 9.59500405313704e-17
45 9.59500405313704e-17
46 9.59500405313704e-17
47 9.59500405313704e-17
48 9.59500405313704e-17
49 9.59500405313704e-17
50 7.46497758975963e-17
};
\addlegendentry{$f \in [0.5 , 1]$ Hz}
\addplot [semithick, black, mark=diamond, mark size=2, mark options={solid}, only marks]
table {%
0 1
1 0.0203338859688381
2 0.00039375935596589
3 3.08031465555279e-05
4 4.53235262419085e-06
5 1.09751598753348e-06
6 3.67014806315636e-07
7 1.50805728475937e-07
8 7.20138010243145e-08
9 3.92866848336345e-08
10 2.36088429850788e-08
11 1.3576300219153e-08
12 8.33102830480686e-09
13 6.83784758601254e-09
14 5.81739107957265e-09
15 3.40538925514246e-09
16 2.65329049419667e-09
17 1.9869913671975e-09
18 1.42306154992201e-09
19 8.41659596685618e-10
20 7.45765600897958e-10
21 7.02898810404598e-10
22 3.98290249128496e-10
23 1.8885960425429e-10
24 1.55724327442289e-10
25 1.11840809410842e-10
26 5.14149328470272e-11
27 2.49807409448052e-11
28 6.15997165958573e-12
29 1.86572407886724e-12
30 1.11036795097701e-12
31 6.91370328087946e-14
32 3.41319148271115e-14
33 2.2363201758684e-14
34 1.00124484058097e-14
35 2.02819079045795e-15
36 1.51574647598997e-15
37 1.02085498179344e-15
38 7.35663516008863e-16
39 6.76208183567884e-16
40 4.73709313453669e-16
41 2.37097342034612e-16
42 1.86486844783258e-16
43 1.64564654492404e-16
44 1.64125663478471e-16
45 1.23074377588517e-16
46 1.22903376153191e-16
47 7.95468386204646e-17
48 7.95468386204646e-17
49 7.95468386204646e-17
50 7.95468386204646e-17
51 7.95468386204646e-17
52 7.95468386204646e-17
53 7.95468386204646e-17
54 7.95468386204646e-17
55 7.95468386204646e-17
56 7.95468386204646e-17
57 7.95468386204646e-17
58 7.95468386204646e-17
59 7.95468386204646e-17
60 7.95468386204646e-17
61 7.95468386204646e-17
62 7.95468386204646e-17
63 7.95468386204646e-17
64 7.95468386204646e-17
65 7.95468386204646e-17
66 7.95468386204646e-17
67 7.95468386204646e-17
68 7.95468386204646e-17
69 7.95468386204646e-17
70 7.95468386204646e-17
71 7.95468386204646e-17
72 7.95468386204646e-17
73 7.95468386204646e-17
74 7.95468386204646e-17
75 7.95468386204646e-17
76 7.95468386204646e-17
77 7.95468386204646e-17
78 7.95468386204646e-17
79 7.95468386204646e-17
};
\addlegendentry{$f \in [1, \; 3]$ Hz}
\addplot [semithick, black, mark=square, mark size=2, mark options={solid}, only marks]
table {%
0 1
1 0.138340173896263
2 0.00376294062532942
3 0.00032260081598158
4 4.62421741501723e-05
5 1.10871706098975e-05
6 3.69604287127299e-06
7 1.51902438037638e-06
8 7.25163302562039e-07
9 3.95399136334086e-07
10 2.37707957858299e-07
11 1.36854267200016e-07
12 8.38514802282797e-08
13 6.88271474705132e-08
14 6.06946042459946e-08
15 3.4355619588786e-08
16 2.67370076358029e-08
17 2.00970900166216e-08
18 1.44505571508203e-08
19 8.64053266236997e-09
20 7.79254260341773e-09
21 7.11638515417939e-09
22 4.97977360051651e-09
23 2.96276531360538e-09
24 1.9320736127674e-09
25 1.76184061249113e-09
26 1.15164087105602e-09
27 7.83919838561793e-10
28 5.3356831969935e-10
29 3.80069963162741e-10
30 2.15630744114883e-10
31 2.27630206362762e-11
32 8.94275945458505e-12
33 1.90779329936048e-12
34 7.82711086957152e-13
35 9.61606254460915e-14
36 2.92216346327567e-14
37 2.32631563843053e-14
38 1.18945561840311e-14
39 2.69853284734875e-15
40 1.37318709132594e-15
41 1.31194111871298e-15
42 1.04296078491935e-15
43 5.97073910345991e-16
44 4.78992951253572e-16
45 3.39519719480296e-16
46 2.34792002907958e-16
47 1.85225046007268e-16
48 1.84123268439647e-16
49 1.4212113741035e-16
50 1.18118558698402e-16
51 1.10911974826039e-16
52 7.23350454377114e-17
53 7.23350454377114e-17
54 7.23350454377114e-17
55 7.23350454377114e-17
56 7.23350454377114e-17
57 7.23350454377114e-17
58 7.23350454377114e-17
59 7.23350454377114e-17
60 7.23350454377114e-17
61 7.23350454377114e-17
62 7.23350454377114e-17
63 7.23350454377114e-17
64 7.23350454377114e-17
65 7.23350454377114e-17
66 7.23350454377114e-17
67 7.23350454377114e-17
68 7.23350454377114e-17
69 7.23350454377114e-17
70 7.23350454377114e-17
71 7.23350454377114e-17
72 7.23350454377114e-17
73 7.23350454377114e-17
74 7.23350454377114e-17
75 7.23350454377114e-17
76 7.23350454377114e-17
77 7.23350454377114e-17
78 7.23350454377114e-17
79 7.23350454377114e-17
};
\addlegendentry{$f \in [3 , \; 5]$ Hz}
\addplot [very thick, darkred, dash pattern=on 11.1pt off 4.8pt, forget plot]
table {%
4 1e-19
4 80
};
\end{axis}

\end{tikzpicture}
		\vspace{-0.1cm}
		\caption{Cantilever beam}
	\end{subfigure}
	\hfill
	\begin{subfigure}[b]{0.49\textwidth}
		% This file was created with tikzplotlib v0.10.1.
\begin{tikzpicture}

\definecolor{darkgray176}{RGB}{176,176,176}
\definecolor{darkred}{RGB}{139,0,0}
\definecolor{darkorange25512714}{RGB}{255,127,14}
\definecolor{forestgreen4416044}{RGB}{44,160,44}
\definecolor{lightgray204}{RGB}{204,204,204}
\definecolor{steelblue31119180}{RGB}{31,119,180}

\begin{axis}[
width=\plotw,
height=\ploth,
legend cell align={left},
legend style={fill opacity=0.8, draw opacity=1, text opacity=1, draw=lightgray204,at={(0.98,0.62)},
	anchor=south east,},
log basis y={10},
tick align=outside,
tick pos=left,
x grid style={darkgray176},
xlabel={r},
xmajorgrids,
xmin=-2, xmax=51,
xtick style={color=black},
y grid style={darkgray176},
ylabel={\(\displaystyle \sigma / \sigma_{max}\)},
ymajorgrids,
ymin=1.41253754462275e-18, ymax=7.07945784384138,
ymode=log,
ytick style={color=black},
ytick={1e-21,1e-18,1e-15,1e-12,1e-09,1e-06,0.001,1,1000,1000000},
%yticklabels={
%  \(\displaystyle {10^{-21}}\),
%  \(\displaystyle {10^{-18}}\),
%  \(\displaystyle {10^{-15}}\),
%  \(\displaystyle {10^{-12}}\),
%  \(\displaystyle {10^{-9}}\),
%  \(\displaystyle {10^{-6}}\),
%  \(\displaystyle {10^{-3}}\),
%  \(\displaystyle {10^{0}}\),
%  \(\displaystyle {10^{3}}\),
%  \(\displaystyle {10^{6}}\)
%},
xticklabel style={font = \small},
yticklabel style={font = \small}
]
\addplot [semithick, black, mark=o, mark size=2, mark options={solid}, only marks]
table[x = r, y = svd1,  col sep=comma] {pics/beam4/csv/ex2_svds.csv};
\addlegendentry{$f \in [50 , 70]$ Hz}
\addplot [semithick, black, mark=diamond, mark size=2, mark options={solid}, only marks]
table[x = r, y = svd2,  col sep=comma] {pics/beam4/csv/ex2_svds.csv};
\addlegendentry{$f \in [140 ,160]$ Hz}
\addplot [semithick, black, mark=square, mark size=2, mark options={solid}, only marks]
table[x = r, y = svd3,  col sep=comma] {pics/beam4/csv/ex2_svds.csv};
\addlegendentry{$f \in [250 , 270]$ Hz}
\addplot [ultra thick, darkred, dash pattern=on 14.8pt off 6.4pt, forget plot]
table {%
3 1e-19
3 51
};
\end{axis}

\end{tikzpicture}
		\vspace{-0.1cm}
		\caption{Overhanging beam}
	\end{subfigure}
\vspace{-0.1cm}
	\caption{Singular values decay of the snapshot matrix $\bX$ that corresponds to the integration of the system \eqref{eq:original_noG} using the input signal $\mathrm{u_m}(t)$ with different frequencies $[f_0, f_1]$\si{\hertz}. The vertical line indicates the reduced order chosen for the model.}
	\label{fig:beam_svds}
\end{figure}
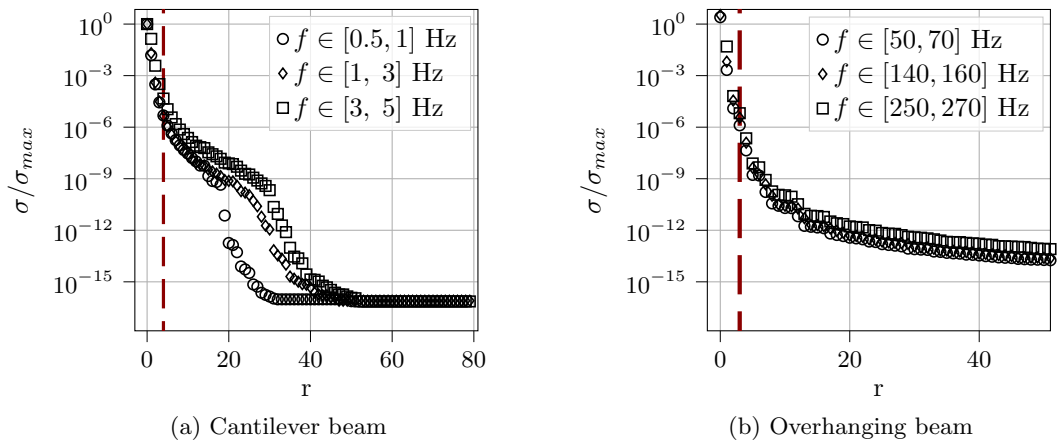 From both plots we can assure that at first there is no significant difference in the singular value decay. Hence, for each beam model we construct three ROMs of the same reduced order r, corresponding to the three frequency ranges, denoted in \Cref{fig:beam_svds}.

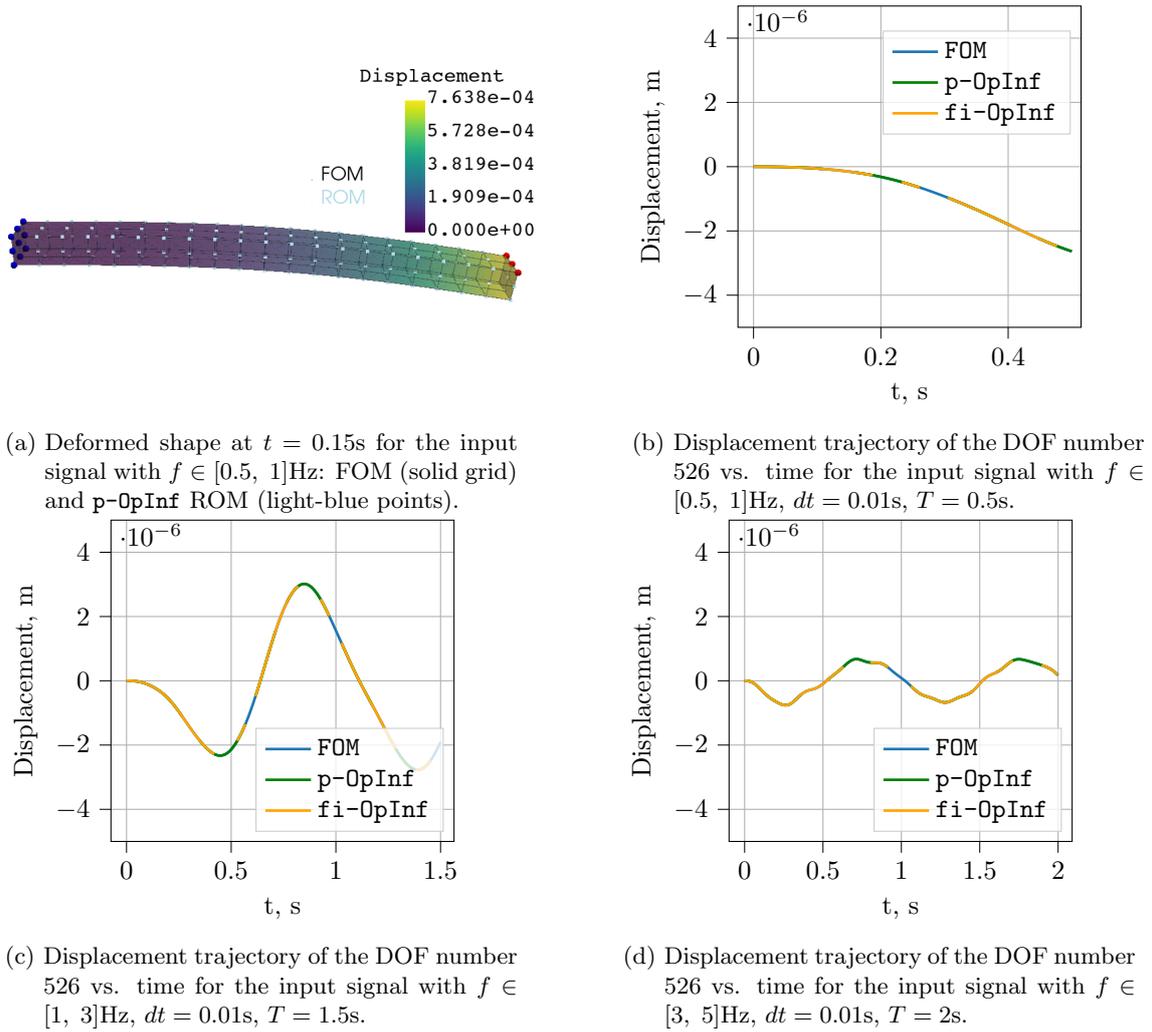
\begin{figure}[h!]
 	\centering
 	\begin{subfigure}[b]{0.45\textwidth}
 		\includegraphics[scale = 0.25]{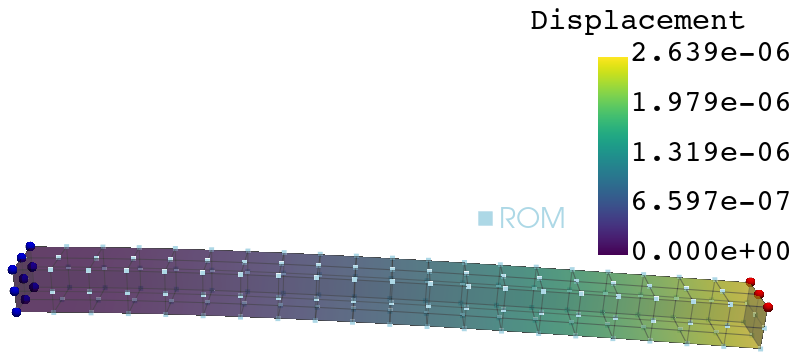}
 		\vspace{1cm}
 		\caption{Deformed shape at $t = 0.15$s for the input signal with $f \in [0.5, \; 1]$\si{\hertz}: FOM (solid grid) and \texttt{p-OpInf} ROM (light-blue points).}
 		\label{fig:ex1_def_shape}
 	\end{subfigure}
 	\hfill
 	\begin{subfigure}[b]{0.45\textwidth}
 		% This file was created with tikzplotlib v0.10.1.
\begin{tikzpicture}

\definecolor{darkgray176}{RGB}{176,176,176}
\definecolor{green}{RGB}{0,128,0}
\definecolor{lightgray204}{RGB}{204,204,204}
\definecolor{orange}{RGB}{255,165,0}
\definecolor{steelblue31119180}{RGB}{31,119,180}

\begin{axis}[
width = \plotw,
height = \ploth,
legend cell align={left},
legend style={
  fill opacity=0.8,
  draw opacity=1,
  text opacity=1,
  at={(0.97,0.6)},
  anchor=south east,
  draw=lightgray204
},
tick align=outside,
tick pos=left,
x grid style={darkgray176},
xlabel={t, s},
xmajorgrids,
xmin=-0.0245, xmax=0.5145,
xtick style={color=black},
y grid style={darkgray176},
ylabel={Displacement, m},
ymajorgrids,
ymin=-5e-06, ymax=5e-06,
ytick style={color=black},
ytick={-6e-06,-4e-06,-2e-06,0,2e-06,4e-06,6e-06},
yticklabels={\ensuremath{-}6,\ensuremath{-}4,\ensuremath{-}2,0,2,4,6},
every y tick scale label/.style={
	at={(-0.07,0.99)},xshift=12pt,anchor=north west,inner sep=0pt
}
]
\addplot [line width=\linew, steelblue31119180]
table[x = time, y = fom,  col sep=comma] {pics/cbeam/csv/ex1_valid_traj_freq_05_1.csv};
\addlegendentry{\texttt{FOM}}
\addplot [line width=\linew, green, dash pattern=on 25.9pt off 11.2pt]
table[x = time, y = popinf,  col sep=comma] {pics/cbeam/csv/ex1_valid_traj_freq_05_1.csv};
\addlegendentry{\texttt{p-OpInf}}
\addplot [line width=\linew, orange, dash pattern=on 44.8pt off 11.2pt on 7pt off 11.2pt]
table[x = time, y = fopinf,  col sep=comma] {pics/cbeam/csv/ex1_valid_traj_freq_05_1.csv};
\addlegendentry{\texttt{fi-OpInf}}
\end{axis}

\end{tikzpicture}
 		\caption{Displacement trajectory of the DOF number 526 vs. time for the input signal with $f \in [0.5, \; 1]$\si{\hertz}, $dt = 0.01$s, $T = 0.5$s.}
 		\label{fig:ex1_valid_0p5_1_traj}
 	\end{subfigure}
 	\vfill
 	\begin{subfigure}[b]{0.45\textwidth}
 		% This file was created with tikzplotlib v0.10.1.
\begin{tikzpicture}

\definecolor{darkgray176}{RGB}{176,176,176}
\definecolor{green}{RGB}{0,128,0}
\definecolor{lightgray204}{RGB}{204,204,204}
\definecolor{orange}{RGB}{255,165,0}
\definecolor{steelblue31119180}{RGB}{31,119,180}

\begin{axis}[
width = \plotw,
height = \ploth,
legend cell align={left},
legend style={
  fill opacity=0.8,
  draw opacity=1,
  text opacity=1,
  at={(0.97,0.03)},
  anchor=south east,
  draw=lightgray204
},
tick align=outside,
tick pos=left,
x grid style={darkgray176},
xlabel={t, s},
xmajorgrids,
xmin=-0.0745, xmax=1.5645,
xtick style={color=black},
y grid style={darkgray176},
ylabel={Displacement, m},
ymajorgrids,
ymin=-5e-06, ymax=5e-06,
ytick style={color=black},
ytick={-6e-06,-4e-06,-2e-06,0,2e-06,4e-06,6e-06},
yticklabels={\ensuremath{-}6,\ensuremath{-}4,\ensuremath{-}2,0,2,4,6},
every y tick scale label/.style={
	at={(-0.07,0.99)},xshift=12pt,anchor=north west,inner sep=0pt
}
]
\addplot [line width=\linew, steelblue31119180]
table[x = time, y = fom,  col sep=comma] {pics/cbeam/csv/ex1_valid_traj_freq_1_3.csv};
\addlegendentry{\texttt{FOM}}
\addplot [line width=\linew, green, dash pattern=on 25.9pt off 11.2pt]
table[x = time, y = popinf,  col sep=comma] {pics/cbeam/csv/ex1_valid_traj_freq_1_3.csv};
\addlegendentry{\texttt{p-OpInf}}
\addplot [line width=\linew, orange, dash pattern=on 44.8pt off 11.2pt on 7pt off 11.2pt]
table[x = time, y = fopinf,  col sep=comma] {pics/cbeam/csv/ex1_valid_traj_freq_1_3.csv};
\addlegendentry{\texttt{fi-OpInf}}
\end{axis}

\end{tikzpicture}
 		\caption{Displacement trajectory of the DOF number 526 vs. time for the input signal with $f \in [1, \; 3]$\si{\hertz}, $dt = 0.01$s, $T = 1.5$s.}
 		\label{fig:ex1_valid_1_3_traj}
 	\end{subfigure}
 	\hfill
 	\begin{subfigure}[b]{0.45\textwidth}
 		% This file was created with tikzplotlib v0.10.1.
\begin{tikzpicture}

\definecolor{darkgray176}{RGB}{176,176,176}
\definecolor{green}{RGB}{0,128,0}
\definecolor{lightgray204}{RGB}{204,204,204}
\definecolor{orange}{RGB}{255,165,0}
\definecolor{steelblue31119180}{RGB}{31,119,180}

\begin{axis}[
width = \plotw,
height = \ploth,
legend cell align={left},
legend style={
  fill opacity=0.8,
  draw opacity=1,
  text opacity=1,
  at={(0.97,0.03)},
  anchor=south east,
  draw=lightgray204
},
tick align=outside,
tick pos=left,
x grid style={darkgray176},
xlabel={t, s},
xmajorgrids,
xmin=-0.0995, xmax=2.0895,
xtick style={color=black},
y grid style={darkgray176},
ylabel={Displacement, m},
ymajorgrids,
ymin=-5e-06, ymax=5e-06,
ytick style={color=black},
ytick={-6e-06,-4e-06,-2e-06,0,2e-06,4e-06,6e-06},
yticklabels={\ensuremath{-}6,\ensuremath{-}4,\ensuremath{-}2,0,2,4,6},
every y tick scale label/.style={
	at={(-0.07,0.99)},xshift=12pt,anchor=north west,inner sep=0pt
}
]
\addplot [line width=\linew, steelblue31119180]
table[x = time, y = fom,  col sep=comma] {pics/cbeam/csv/ex1_valid_traj_freq_3_5.csv};
\addlegendentry{\texttt{FOM}}
\addplot [line width=\linew, green, dash pattern=on 25.9pt off 11.2pt]
table[x = time, y = popinf,  col sep=comma] {pics/cbeam/csv/ex1_valid_traj_freq_3_5.csv};
\addlegendentry{\texttt{p-OpInf}}
\addplot [line width=\linew, orange, dash pattern=on 44.8pt off 11.2pt on 7pt off 11.2pt]
table[x = time, y = fopinf,  col sep=comma] {pics/cbeam/csv/ex1_valid_traj_freq_3_5.csv};
\addlegendentry{\texttt{fi-OpInf}}
\end{axis}

\end{tikzpicture}
 		\caption{Displacement trajectory of the DOF number 526 vs. time for the input signal with $f \in [3, \; 5]$\si{\hertz}, $dt = 0.01$s, $T = 2$s.}
 		\label{fig:ex1_valid_3_5traj}
 	\end{subfigure}
 	\vspace{-0.2cm}
 	\caption{Validation of the ROM $r = 4$ results for the cantilever beam model.}
 	\label{fig:cbeam_validation}
 \end{figure}

 In \Cref{fig:ex1_def_shape} the deformed shape for the FOM and \texttt{p-OpInf}-ROM is shown to prove the meaningful 3D behavior of the ROM. In \Cref{fig:ex1_valid_0p5_1_traj,fig:ex1_valid_1_3_traj,fig:ex1_valid_3_5traj}, the trajectories of one DOF for different input signals $\mathrm{u_m}$ are plotted for FOM, \texttt{p-OpInf}-ROM, and \texttt{fi-OpInf}-ROM. The figures show consistency and sufficient accuracy of results. The respective error plots are given in \Cref{fig:cbeam_valid_errors}. 
 
\begin{figure}[h!]
	\centering
	\begin{subfigure}[b]{0.3\textwidth}
		% This file was created with tikzplotlib v0.10.1.
\begin{tikzpicture}

\definecolor{darkgray176}{RGB}{176,176,176}
\definecolor{darkorange}{RGB}{255,140,0}
\definecolor{green01270}{RGB}{0,127,0}
\definecolor{lightgray204}{RGB}{204,204,204}

\begin{axis}[
width=\plotwtripple,
height=\plothtripple,
legend cell align={left},
legend style={
  fill opacity=0.8,
  draw opacity=1,
  text opacity=1,
  at={(2.4,1.06)},
  anchor=south east,
  draw=lightgray204,
  legend columns=2,
  nodes={inner xsep=20pt}
},
log basis y={10},
tick align=outside,
tick pos=left,
x grid style={darkgray176},
xlabel={t, s},
xmajorgrids,
xmin=-0.025, xmax=0.525,
xtick style={color=black},
y grid style={darkgray176},
ylabel={Relative error $\epsilon_{\mathrm{err}}(t)$},
ylabel style = {yshift = -7pt},
ymajorgrids,
ymin=1e-07, ymax=0.01,
ymode=log,
ytick style={color=black, font=\small},
ytick={1e-08,1e-07,1e-06,1e-05,0.0001,0.001,0.01,0.1,1},
yticklabels={
  \(\displaystyle {10^{-8}}\),
  \(\displaystyle {10^{-7}}\),
  \(\displaystyle {10^{-6}}\),
  \(\displaystyle {10^{-5}}\),
  \(\displaystyle {10^{-4}}\),
  \(\displaystyle {10^{-3}}\),
  \(\displaystyle {10^{-2}}\),
  \(\displaystyle {10^{-1}}\),
  \(\displaystyle {10^{0}}\)
},
xticklabel style={font = \small},
yticklabel style={font = \small}
]
\addplot [line width=\linew, green01270]
table[x = time, y expr={\thisrow{popinf} <= \pgfkeysvalueof{/pgfplots/ymin} ? \pgfkeysvalueof{/pgfplots/ymin} : \thisrow{popinf}},col sep=comma] {pics/cbeam/csv/ex1_valid_error_freq_05_1.csv};
\addlegendentry{\texttt{p-OpInf}}
\addplot [line width=\linew, darkorange]
table[x = time, y expr={\thisrow{fopinf} <= \pgfkeysvalueof{/pgfplots/ymin} ? \pgfkeysvalueof{/pgfplots/ymin} : \thisrow{fopinf}}, col sep=comma] {pics/cbeam/csv/ex1_valid_error_freq_05_1.csv};
\addlegendentry{\texttt{fi-OpInf}}
\node at (axis cs:0.2,0.0005)  [draw, rectangle, fill = white, anchor=south] {$f \in [0.5, \; 1]$ Hz};
\end{axis}

\end{tikzpicture}
	\end{subfigure}
	\hspace{0.7cm}
	\begin{subfigure}[b]{0.3\textwidth}
		% This file was created with tikzplotlib v0.10.1.
\begin{tikzpicture}

\definecolor{darkgray176}{RGB}{176,176,176}
\definecolor{darkorange}{RGB}{255,140,0}
\definecolor{green01270}{RGB}{0,127,0}
\definecolor{lightgray204}{RGB}{204,204,204}

\begin{axis}[
width=\plotwtripple,
height=\plothtripple,
%legend cell align={left},
%legend style={
%  fill opacity=0.8,
%  draw opacity=1,
%  text opacity=1,
%  at={(0.97,0.03)},
%  anchor=south east,
%  draw=lightgray204
%},
log basis y={10},
tick align=outside,
tick pos=left,
x grid style={darkgray176},
xlabel={t, s},
xmajorgrids,
xmin=-0.075, xmax=1.575,
xtick style={color=black,font=\small},
y grid style={darkgray176},
ymajorgrids,
ymin=1e-07, ymax=0.01,
ymode=log,
ytick style={color=black},
yticklabel style = {font=\small,xshift = 7pt},
ytick={1e-08,1e-07,1e-06,1e-05,0.0001,0.001,0.01,0.1,1},
yticklabels={
  \(\displaystyle {10^{-8}}\),
  \(\displaystyle {10^{-7}}\),
  \(\displaystyle {10^{-6}}\),
  \(\displaystyle {10^{-5}}\),
  \(\displaystyle {10^{-4}}\),
  \(\displaystyle {10^{-3}}\),
  \(\displaystyle {10^{-2}}\),
  \(\displaystyle {10^{-1}}\),
  \(\displaystyle {10^{0}}\)
},
xticklabel style={font = \small},
yticklabel style={font = \small}
]
\addplot [line width=\linew, each nth point=2,green01270]
table[x = time, y expr={\thisrow{popinf} <= \pgfkeysvalueof{/pgfplots/ymin} ? \pgfkeysvalueof{/pgfplots/ymin} : \thisrow{popinf}},  col sep=comma] {pics/cbeam/csv/ex1_valid_error_freq_1_3.csv};
%\addlegendentry{\texttt{p-OpInf}}
\addplot [line width=\linew,each nth point=2, darkorange]
table[x = time, y expr={\thisrow{fopinf} <= \pgfkeysvalueof{/pgfplots/ymin} ? \pgfkeysvalueof{/pgfplots/ymin} : \thisrow{fopinf}},  col sep=comma] {pics/cbeam/csv/ex1_valid_error_freq_1_3.csv};
%\addlegendentry{\texttt{fi-OpInf}}
\node at (axis cs:0.7,0.0005)  [draw, rectangle, fill = white, anchor=south] {$f \in [1, \; 3]$ Hz};
\end{axis}

\end{tikzpicture}
	\end{subfigure}
	\hfill
	\begin{subfigure}[b]{0.3\textwidth}
		% This file was created with tikzplotlib v0.10.1.
\begin{tikzpicture}

\definecolor{darkgray176}{RGB}{176,176,176}
\definecolor{darkorange}{RGB}{255,140,0}
\definecolor{green01270}{RGB}{0,127,0}
\definecolor{lightgray204}{RGB}{204,204,204}

\begin{axis}[
width=\plotwtripple,
height=\plothtripple,
%legend cell align={left},
%legend style={
%  fill opacity=0.8,
%  draw opacity=1,
%  text opacity=1,
%  at={(0.97,0.03)},
%  anchor=south east,
%  draw=lightgray204
%},
log basis y={10},
tick align=outside,
tick pos=left,
x grid style={darkgray176},
xlabel={t, s},
xmajorgrids,
xmin=-0.1, xmax=2.1,
xtick style={color=black},
xtick={-0.25,0,0.25,0.5,0.75,1,1.25,1.5,1.75,2,2.25},
xticklabels={-0.2,0.0,0.2,0.5,0.8,1.0,1.2,1.5,1.8,2.0,2.2},
y grid style={darkgray176},
ymajorgrids,
ymin=1e-07, ymax=0.01,
ymode=log,
ytick style={color=black},
yticklabel style = {font = \small, xshift = 7pt},
xticklabel style = {font = \small},
xtick={0,0.5,1.0,1.5,2},
xticklabels = {0, 0.5, 1.0, 1.5, 2.0},
ytick={1e-08,1e-07,1e-06,1e-05,0.0001,0.001,0.01,0.1,1},
yticklabels={
  \(\displaystyle {10^{-8}}\),
  \(\displaystyle {10^{-7}}\),
  \(\displaystyle {10^{-6}}\),
  \(\displaystyle {10^{-5}}\),
  \(\displaystyle {10^{-4}}\),
  \(\displaystyle {10^{-3}}\),
  \(\displaystyle {10^{-2}}\),
  \(\displaystyle {10^{-1}}\),
  \(\displaystyle {10^{0}}\)
},
xticklabel style={font = \small},
yticklabel style={font = \small}
]
\addplot [line width=\linew,each nth point=2, green01270]
table[x = time, y expr={\thisrow{popinf} <= \pgfkeysvalueof{/pgfplots/ymin} ? \pgfkeysvalueof{/pgfplots/ymin} : \thisrow{popinf}},  col sep=comma] {pics/cbeam/csv/ex1_valid_error_freq_3_5.csv};
%\addlegendentry{\texttt{p-OpInf}}
\addplot [line width=\linew,each nth point=2, darkorange]
table[x = time, y expr={\thisrow{fopinf} <= \pgfkeysvalueof{/pgfplots/ymin} ? \pgfkeysvalueof{/pgfplots/ymin} : \thisrow{fopinf}},  col sep=comma] {pics/cbeam/csv/ex1_valid_error_freq_3_5.csv};
%\addlegendentry{\texttt{fi-OpInf}}
\node at (axis cs:1,0.0000005)  [draw, rectangle, fill = white, anchor=south] {$f \in [3, \; 5]$ Hz};
\end{axis}

\end{tikzpicture}
	\end{subfigure}
\vspace{-0.75cm}
	\caption{Cantilever beam: Relative error vs. time for ROMs $r = 4$ constructed with different frequency ranges $[f_0, f_1 ]$.}
	\label{fig:cbeam_valid_errors}
\end{figure} \FloatBarrier

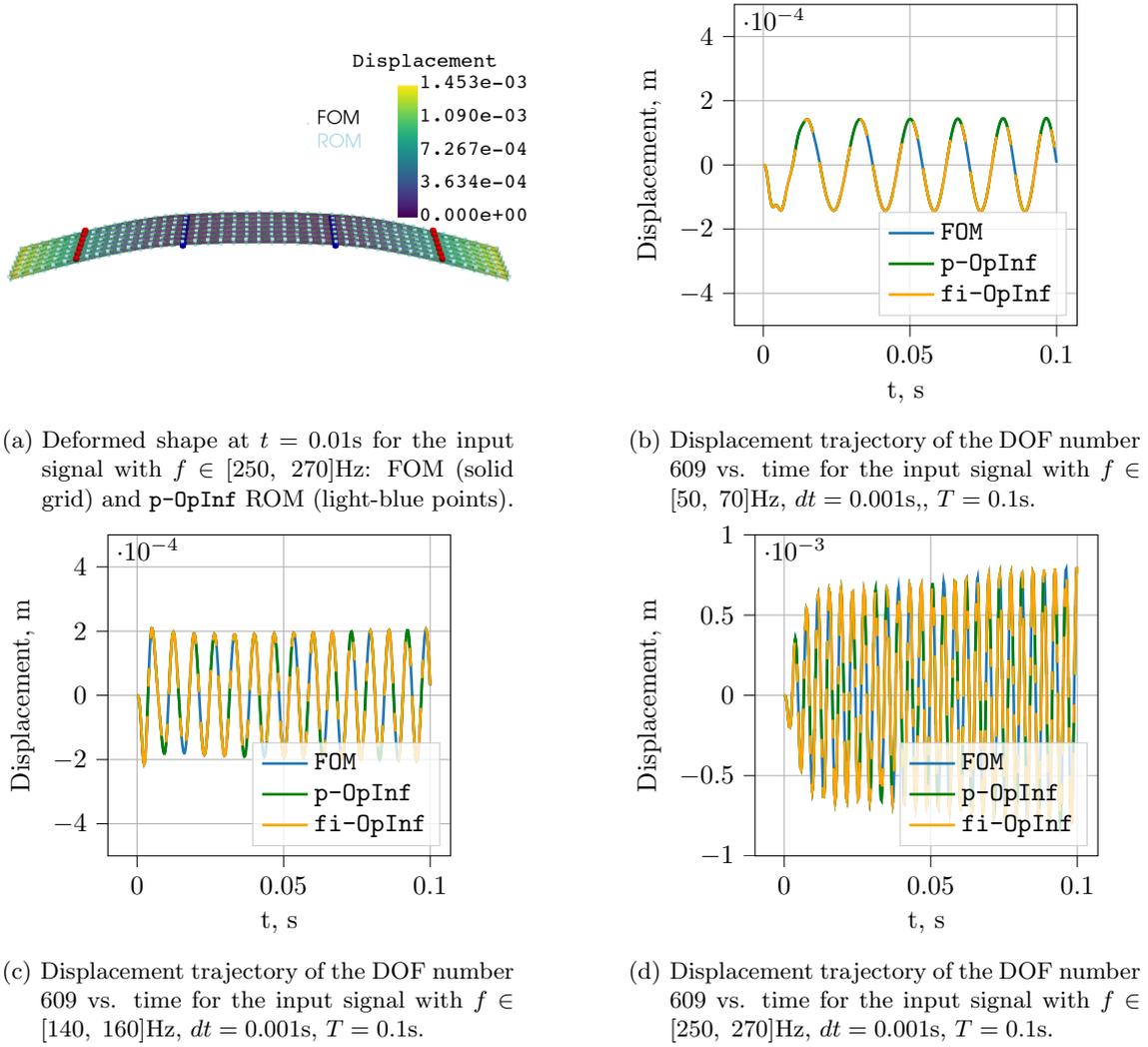
\begin{figure}[h!]
	\centering
	\begin{subfigure}[b]{0.45\textwidth}
			\vspace*{0.5cm}
		\includegraphics[scale = 0.25]{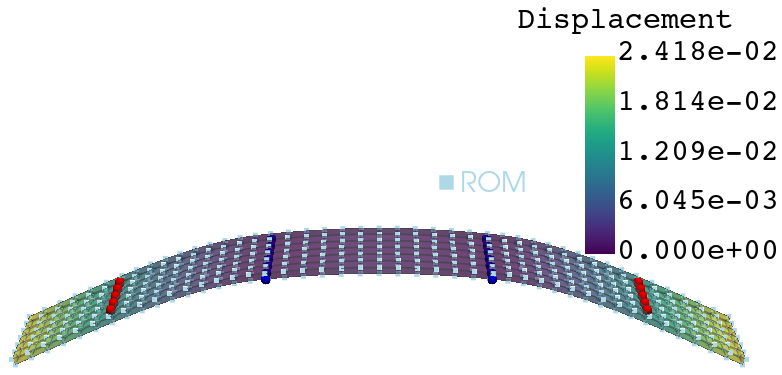}
		\vspace{1.2cm}
		\caption{Deformed shape at $t = 0.01$s for the input signal with $f \in [250, \; 270]$\si{\hertz}: FOM (solid grid) and \texttt{p-OpInf} ROM (light-blue points).}
		\label{fig:ex2_def_shape}
	\end{subfigure}
	\hfill
	\begin{subfigure}[b]{0.45\textwidth}
		% This file was created with tikzplotlib v0.10.1.
\begin{tikzpicture}

\definecolor{darkgray176}{RGB}{176,176,176}
\definecolor{green}{RGB}{0,128,0}
\definecolor{lightgray204}{RGB}{204,204,204}
\definecolor{orange}{RGB}{255,165,0}
\definecolor{steelblue31119180}{RGB}{31,119,180}

\begin{axis}[
width = \plotw,
height = \ploth,
legend cell align={left},
legend style={
  fill opacity=0.8,
  draw opacity=1,
  text opacity=1,
  at={(0.97,0.03)},
  anchor=south east,
  draw=lightgray204
},
tick align=outside,
tick pos=left,
x grid style={darkgray176},
xlabel={t, s},
xmajorgrids,
xmin=-0.009875, xmax=0.107,
xtick style={color=black},
y grid style={darkgray176},
ylabel={Displacement, m},
ymajorgrids,
ymin=-0.0005, ymax=0.0005,
ytick style={color=black},
scaled x ticks=false,                          % disables sci notation on x-axis
x tick label style={/pgf/number format/fixed}, % formats x-axis ticks
every y tick scale label/.style={
	at={(-0.07,0.99)},xshift=12pt,anchor=north west,inner sep=0pt
}
]
\addplot [line width=\linew, steelblue31119180]
table[x = time, y = fom,  col sep=comma] {pics/beam4/csv/ex2_valid_traj_freq_50_70.csv};
\addlegendentry{\texttt{FOM}}
\addplot [line width=\linew, green, dash pattern=on 25.9pt off 11.2pt]
table[x = time, y = popinf,  col sep=comma] {pics/beam4/csv/ex2_valid_traj_freq_50_70.csv};
\addlegendentry{\texttt{p-OpInf}}
\addplot [line width=\linew, orange, dash pattern=on 44.8pt off 11.2pt on 7pt off 11.2pt]
table[x = time, y = fopinf,  col sep=comma] {pics/beam4/csv/ex2_valid_traj_freq_50_70.csv};
\addlegendentry{\texttt{fi-OpInf}}
\end{axis}

\end{tikzpicture}
		\caption{Displacement trajectory of the DOF number 609 vs. time for the input signal with $f \in [50, \; 70]$\si{\hertz}, $dt = 0.001$s,, $T = 0.1$s.}
		\label{fig:ex2_valid_0p5_1_traj}
	\end{subfigure}
	\vfill
	\begin{subfigure}[b]{0.45\textwidth}
		% This file was created with tikzplotlib v0.10.1.
\begin{tikzpicture}

\definecolor{darkgray176}{RGB}{176,176,176}
\definecolor{green}{RGB}{0,128,0}
\definecolor{lightgray204}{RGB}{204,204,204}
\definecolor{orange}{RGB}{255,165,0}
\definecolor{steelblue31119180}{RGB}{31,119,180}

\begin{axis}[
width = \plotw,
height = \ploth,
legend cell align={left},
legend style={
  fill opacity=0.8,
  draw opacity=1,
  text opacity=1,
  at={(0.97,0.03)},
  anchor=south east,
  draw=lightgray204
},
tick align=outside,
tick pos=left,
x grid style={darkgray176},
scaled x ticks=false,                          % disables sci notation on x-axis
x tick label style={/pgf/number format/fixed}, % formats x-axis ticks
xlabel={t, s},
xmajorgrids,
xmin=-0.00995, xmax=0.107,
xtick style={color=black},
y grid style={darkgray176},
ylabel={Displacement, m},
ymajorgrids,
ymin=-0.0005, ymax=0.0005,
ytick style={color=black},
every y tick scale label/.style={
	at={(-0.07,0.99)},xshift=12pt,anchor=north west,inner sep=0pt
}
]
\addplot [line width=\linew, steelblue31119180]
table[x = time, y = fom,  col sep=comma] {pics/beam4/csv/ex2_valid_traj_freq_140_160.csv};
\addlegendentry{\texttt{FOM}}
\addplot [line width=\linew, green, dash pattern=on 25.9pt off 11.2pt]
table[x = time, y = popinf,  col sep=comma] {pics/beam4/csv/ex2_valid_traj_freq_140_160.csv};
\addlegendentry{\texttt{p-OpInf}}
\addplot [line width=\linew, orange, dash pattern=on 44.8pt off 11.2pt on 7pt off 11.2pt]
table[x = time, y = fopinf,  col sep=comma] {pics/beam4/csv/ex2_valid_traj_freq_140_160.csv};
\addlegendentry{\texttt{fi-OpInf}}
\end{axis}

\end{tikzpicture}
		\caption{Displacement trajectory of the DOF number 609 vs. time for the input signal with $f \in [140, \; 160]$\si{\hertz}, $dt = 0.001$s, $T = 0.1$s.}
		\label{fig:ex2_valid_1_3_traj}
	\end{subfigure}
	\hfill
	\begin{subfigure}[b]{0.45\textwidth}
		% This file was created with tikzplotlib v0.10.1.
\begin{tikzpicture}

\definecolor{darkgray176}{RGB}{176,176,176}
\definecolor{green}{RGB}{0,128,0}
\definecolor{lightgray204}{RGB}{204,204,204}
\definecolor{orange}{RGB}{255,165,0}
\definecolor{steelblue31119180}{RGB}{31,119,180}

\begin{axis}[
width = \plotw,
height = \ploth,
legend cell align={left},
legend style={
  fill opacity=0.8,
  draw opacity=1,
  text opacity=1,
  at={(0.97,0.03)},
  anchor=south east,
  draw=lightgray204
},
tick align=outside,
tick pos=left,
x grid style={darkgray176},
xlabel={t, s},
xmajorgrids,
xmin=-0.00995, xmax=0.107,
xtick style={color=black},
scaled x ticks=false,                          % disables sci notation on x-axis
x tick label style={/pgf/number format/fixed}, % formats x-axis ticks
y grid style={darkgray176},
ylabel={Displacement, m},
ymajorgrids,
ymin=-0.001, ymax=0.001,
ytick style={color=black},
every y tick scale label/.style={
	at={(-0.07,0.99)},xshift=12pt,anchor=north west,inner sep=0pt
}
]
\addplot [line width=\linew, steelblue31119180]
table[x = time, y = fom,  col sep=comma] {pics/beam4/csv/ex2_valid_traj_freq_250_270.csv};
\addlegendentry{\texttt{FOM}}
\addplot [line width=\linew, green, dash pattern=on 25.9pt off 11.2pt]
table[x = time, y = popinf,  col sep=comma] {pics/beam4/csv/ex2_valid_traj_freq_250_270.csv};
\addlegendentry{\texttt{p-OpInf}}
\addplot [line width=\linew, orange, dash pattern=on 44.8pt off 11.2pt on 7pt off 11.2pt]
table[x = time, y = fopinf,  col sep=comma] {pics/beam4/csv/ex2_valid_traj_freq_250_270.csv};
\addlegendentry{\texttt{fi-OpInf}}
\end{axis}

\end{tikzpicture}
		\caption{Displacement trajectory of the DOF number 609 vs. time for the input signal with $f \in [250, \; 270]$\si{\hertz}, $dt = 0.001$s, $T = 0.1$s.}
		\label{fig:ex2_valid_3_5traj}
	\end{subfigure}
	\caption{Validation of the ROM $r = 3$ results for the overhanging beam model.}
	\label{fig:beam4_validation}
\end{figure}

\begin{figure}[h!]
	\centering
	\begin{subfigure}[b]{0.3\textwidth}
		% This file was created with tikzplotlib v0.10.1.
\begin{tikzpicture}

\definecolor{darkgray176}{RGB}{176,176,176}
\definecolor{darkorange}{RGB}{255,140,0}
\definecolor{green01270}{RGB}{0,127,0}
\definecolor{lightgray204}{RGB}{204,204,204}

\begin{axis}[
width=\plotwtripple,
height=\plothtripple,
legend cell align={left},
legend style={
	fill opacity=0.8,
	draw opacity=1,
	text opacity=1,
	at={(2.4,1.06)},
	anchor=south east,
	draw=lightgray204,
	legend columns=2,
	nodes={inner xsep=20pt} 
},
log basis y={10},
tick align=outside,
tick pos=left,
x grid style={darkgray176},
xlabel={t (s)},
xmajorgrids,
xmin=-0.01, xmax=0.11,
xtick style={color=black},
xtick={0,0.05,0.1,0.15,0.2,0.225},
xticklabels={0,0.05,0.1,0.15,0.2,0.2},
y grid style={darkgray176},
ylabel={Relative error $\epsilon_{\mathrm{err}}(t)$},
ylabel style = {yshift = -7pt},
ymajorgrids,
ymin=1e-08, ymax=0.01,
ymode=log,
ytick style={color=black},
yticklabel style={color=black, font=\small, xshift = 6pt},
ytick={1e-08,1e-07,1e-06,1e-05,0.0001,0.001,0.01,0.1,1},
yticklabels={
  \(\displaystyle {10^{-8}}\),
  \(\displaystyle {10^{-7}}\),
  \(\displaystyle {10^{-6}}\),
  \(\displaystyle {10^{-5}}\),
  \(\displaystyle {10^{-4}}\),
  \(\displaystyle {10^{-3}}\),
  \(\displaystyle {10^{-2}}\),
  \(\displaystyle {10^{-1}}\),
  \(\displaystyle {10^{0}}\)
}
]
\addplot [line width=\linew, each nth point=2, green01270]
table[x = time, y expr={\thisrow{popinf} <= \pgfkeysvalueof{/pgfplots/ymin} ? \pgfkeysvalueof{/pgfplots/ymin} : \thisrow{popinf}},  col sep=comma] {pics/beam4/csv/ex2_valid_error_freq_50_70.csv};
\addlegendentry{\texttt{p-OpInf}}
\addplot [line width=\linew, each nth point=2, darkorange]
table[x = time, y expr={\thisrow{fopinf} <= \pgfkeysvalueof{/pgfplots/ymin} ? \pgfkeysvalueof{/pgfplots/ymin} : \thisrow{fopinf}},  col sep=comma] {pics/beam4/csv/ex2_valid_error_freq_50_70.csv};
\addlegendentry{\texttt{fi-OpInf}}
\node at (axis cs:0.06,0.0005)  [draw, rectangle, fill = white, anchor=south] {$f \in [50, \; 70]$ Hz};
\end{axis}

\end{tikzpicture}
	\end{subfigure}
	\hspace{0.7cm}
	\begin{subfigure}[b]{0.3\textwidth}
		% This file was created with tikzplotlib v0.10.1.
\begin{tikzpicture}

\definecolor{darkgray176}{RGB}{176,176,176}
\definecolor{darkorange}{RGB}{255,140,0}
\definecolor{green01270}{RGB}{0,127,0}
\definecolor{lightgray204}{RGB}{204,204,204}

\begin{axis}[
width=\plotwtripple,
height=\plothtripple,
%legend cell align={left},
%legend style={
%  fill opacity=0.8,
%  draw opacity=1,
%  text opacity=1,
%  at={(0.97,0.03)},
%  anchor=south east,
%  draw=lightgray204
%},
log basis y={10},
tick align=outside,
tick pos=left,
x grid style={darkgray176},
xlabel={t, s},
xmajorgrids,
xmin=-0.01, xmax=0.11,
xtick style={color=black},
xtick={0,0.05,0.1,0.15,0.2,0.225},
xticklabels={0,0.05,0.1,0.15,0.2,0.2},
y grid style={darkgray176},
ymajorgrids,
ymin=1e-08, ymax=0.01,
ymode=log,
ytick style={color=black},
yticklabel style={color=black, font=\small, xshift = 6pt},
ytick={1e-08,1e-07,1e-06,1e-05,0.0001,0.001,0.01,0.1,1},
yticklabels={
  \(\displaystyle {10^{-8}}\),
  \(\displaystyle {10^{-7}}\),
  \(\displaystyle {10^{-6}}\),
  \(\displaystyle {10^{-5}}\),
  \(\displaystyle {10^{-4}}\),
  \(\displaystyle {10^{-3}}\),
  \(\displaystyle {10^{-2}}\),
  \(\displaystyle {10^{-1}}\),
  \(\displaystyle {10^{0}}\)
}
]
\addplot [line width=\linew, each nth point=2, green01270]
table[x = time, y expr={\thisrow{popinf} <= \pgfkeysvalueof{/pgfplots/ymin} ? \pgfkeysvalueof{/pgfplots/ymin} : \thisrow{popinf}},  col sep=comma] {pics/beam4/csv/ex2_valid_error_freq_140_160.csv};
%\addlegendentry{p-OpInf}
\addplot [line width=\linew, each nth point=2, darkorange]
table[x = time, y expr={\thisrow{fopinf} <= \pgfkeysvalueof{/pgfplots/ymin} ? \pgfkeysvalueof{/pgfplots/ymin} : \thisrow{fopinf}},  col sep=comma] {pics/beam4/csv/ex2_valid_error_freq_140_160.csv};
%\addlegendentry{fi-OpInf}
\node at (axis cs:0.06,0.0007)  [draw, rectangle, fill = white, anchor=south] {$f \in [140, \; 160]$ Hz};
\end{axis}

\end{tikzpicture}
	\end{subfigure}
	\hfill
	\begin{subfigure}[b]{0.3\textwidth}
		% This file was created with tikzplotlib v0.10.1.
\begin{tikzpicture}

\definecolor{darkgray176}{RGB}{176,176,176}
\definecolor{darkorange}{RGB}{255,140,0}
\definecolor{green01270}{RGB}{0,127,0}
\definecolor{lightgray204}{RGB}{204,204,204}

\begin{axis}[
width=\plotwtripple,
height=\plothtripple,
%legend cell align={left},
%legend style={
%  fill opacity=0.8,
%  draw opacity=1,
%  text opacity=1,
%  at={(0.97,0.03)},
%  anchor=south east,
%  draw=lightgray204
%},
log basis y={10},
tick align=outside,
tick pos=left,
x grid style={darkgray176},
xlabel={t, s},
xmajorgrids,
xmin=-0.01, xmax=0.11,
xtick style={color=black},
xtick={0,0.05,0.1,0.15,0.2,0.225},
xticklabels={0,0.05,0.1,0.15,0.2,0.2},
y grid style={darkgray176},
ymajorgrids,
ymin=1e-08, ymax=0.01,
ymode=log,
ytick style={color=black},
yticklabel style={color=black, font=\small, xshift = 6pt},
ytick={1e-08,1e-07,1e-06,1e-05,0.0001,0.001,0.01,0.1,1},
yticklabels={
  \(\displaystyle {10^{-8}}\),
  \(\displaystyle {10^{-7}}\),
  \(\displaystyle {10^{-6}}\),
  \(\displaystyle {10^{-5}}\),
  \(\displaystyle {10^{-4}}\),
  \(\displaystyle {10^{-3}}\),
  \(\displaystyle {10^{-2}}\),
  \(\displaystyle {10^{-1}}\),
  \(\displaystyle {10^{0}}\)
}
]
\addplot [line width=\linew, each nth point=2, green01270]
table[x = time,y expr={\thisrow{popinf} <= \pgfkeysvalueof{/pgfplots/ymin} ? \pgfkeysvalueof{/pgfplots/ymin} : \thisrow{popinf}},  col sep=comma] {pics/beam4/csv/ex2_valid_error_freq_250_270.csv};
%\addlegendentry{p-OpInf}
\addplot [line width=\linew, each nth point=2, darkorange]
table[x = time, y expr={\thisrow{fopinf} <= \pgfkeysvalueof{/pgfplots/ymin} ? \pgfkeysvalueof{/pgfplots/ymin} : \thisrow{fopinf}},  col sep=comma] {pics/beam4/csv/ex2_valid_error_freq_250_270.csv};
%\addlegendentry{fi-OpInf}
\node at (axis cs:0.06,0.001)  [draw, rectangle, fill = white, anchor=south] {$f \in [250, \; 270]$ Hz};
\end{axis}

\end{tikzpicture}
	\end{subfigure}
\vspace{-0.7cm}
	\caption{Overhanging beam: Relative error vs. time for ROMs $r = 3$ constructed with different frequency ranges $[f_0, f_1 ]$.}
	\label{fig:beam4_valid_errors}
\end{figure} \FloatBarrier

The validation plots for the overhanging beam are presented in \Cref{fig:beam4_validation} and \Cref{fig:beam4_valid_errors}.
To get a rough understanding, how the input signals $\mathrm{u_m}$ with different frequencies for the reduced modeling affect the results quality, we performed the test simulation under another loading condition. More precisely, for each beam model we test the resulting ROMs applying a \textit{mono-frequent} input signal $\mathrm{u}(t)$ with a frequency inside and outside the interval $[f_0 , f_1] \si{\hertz}$.

 We test the cantilever beam for $t \in [0 , 16]s$ with the time step $dt = 0.001s$ at different frequencies, as shown on \Cref{fig:cbeam_tend2_r4}. The black lines denote the frequencies for the chirp signal $\mathrm{u_m}$, used for generating snapshots and constructing the reduced-order model.

\begin{figure}[h!]
	\centering
	\begin{subfigure}[b]{0.3\textwidth}
		% This file was created with tikzplotlib v0.10.1.
\begin{tikzpicture}

\definecolor{darkgray176}{RGB}{176,176,176}
\definecolor{darkorange}{RGB}{255,140,0}
\definecolor{green01270}{RGB}{0,127,0}
\definecolor{lightgray204}{RGB}{204,204,204}

\begin{axis}[
width=\plotwtripple,
height=\plothtripple,
legend cell align={left},
legend style={
  fill opacity=0.8,
  draw opacity=1,
  text opacity=1,
  at={(0.97,0.03)},
  anchor=south east,
  draw=lightgray204
},
log basis y={10},
tick align=outside,
tick pos=left,
x grid style={darkgray176},
xlabel={$f$, Hz},
xmajorgrids,
xmin=-0.7395, xmax=15.7495,
xtick style={color=black},
y grid style={darkgray176},
ylabel={Max. relative error},
ymajorgrids,
ymin=1e-06, 
ymax=0.5,
xmax = 8,
ymode=log,
yticklabel style = {font = \small, xshift = 7pt},
xticklabel style = {font = \small},
ytick style={color=black},
ytick={1e-07,1e-06,1e-05,0.0001,0.001,0.01,0.1,1,10},
yticklabels={
  \(\displaystyle {10^{-7}}\),
  \(\displaystyle {10^{-6}}\),
  \(\displaystyle {10^{-5}}\),
  \(\displaystyle {10^{-4}}\),
  \(\displaystyle {10^{-3}}\),
  \(\displaystyle {10^{-2}}\),
  \(\displaystyle {10^{-1}}\),
  \(\displaystyle {10^{0}}\),
  \(\displaystyle {10^{1}}\)
}
]
\addplot[green01270, mark=triangle*, mark size=2.5, mark options={solid}, only marks] table[x = fr, y = popinf,  col sep=comma] {pics/cbeam/csv/ex1_text_freq_05_1_popinf.csv};
\addlegendentry{\texttt{p-OpInf}}
\addplot [semithick, darkorange, mark=*, mark size=2.5, mark options={solid}, only marks]
table[x = fr, y = fopinf,  col sep=comma] {pics/cbeam/csv/ex1_text_freq_05_1_fopinf.csv};
\addlegendentry{\texttt{fi-OpInf}}
\addplot [ultra thick, black, forget plot]
table {%
0.5 1e-06
0.5 1
};
\addplot [ultra thick, black, forget plot]
table {%
	1 1e-06
	1 1
};
%\legend{};
\end{axis}

\end{tikzpicture}
	\end{subfigure}
\hspace{0.4cm}
	\begin{subfigure}[b]{0.3\textwidth}
		% This file was created with tikzplotlib v0.10.1.
\begin{tikzpicture}

\definecolor{darkgray176}{RGB}{176,176,176}
\definecolor{darkorange}{RGB}{255,140,0}
\definecolor{green01270}{RGB}{0,127,0}
\definecolor{lightgray204}{RGB}{204,204,204}

\begin{axis}[
width=\plotwtripple,
height=\plothtripple,
legend cell align={left},
legend style={
  fill opacity=0.8,
  draw opacity=1,
  text opacity=1,
  at={(0.97,0.03)},
  anchor=south east,
  draw=lightgray204
},
log basis y={10},
tick align=outside,
tick pos=left,
x grid style={darkgray176},
xlabel={$f$, Hz},
xmajorgrids,
xmin=-0.7395, xmax=15.7495,
xtick style={color=black},
y grid style={darkgray176},
ylabel={\empty},
ymajorgrids,
ymin=1e-06, 
ymax=0.5,
xmax = 8,
yticklabel style = {font = \small, xshift = 7pt},
xticklabel style = {font = \small},
ymode=log,
ytick style={color=black},
ytick={1e-07,1e-06,1e-05,0.0001,0.001,0.01,0.1,1,10,100},
yticklabels={
  \(\displaystyle {10^{-7}}\),
  \(\displaystyle {10^{-6}}\),
  \(\displaystyle {10^{-5}}\),
  \(\displaystyle {10^{-4}}\),
  \(\displaystyle {10^{-3}}\),
  \(\displaystyle {10^{-2}}\),
  \(\displaystyle {10^{-1}}\),
  \(\displaystyle {10^{0}}\),
  \(\displaystyle {10^{1}}\),
  \(\displaystyle {10^{2}}\)
}
]
\addplot [semithick, green01270, mark=triangle*, mark size=2.5, mark options={solid}, only marks]
table[x = fr, y = popinf,  col sep=comma] {pics/cbeam/csv/ex1_text_freq_1_3_popinf.csv};
\addlegendentry{p-OpInf}
\addplot [semithick, darkorange, mark=*, mark size=2.5, mark options={solid}, only marks]
table[x = fr, y = fopinf,  col sep=comma] {pics/cbeam/csv/ex1_text_freq_1_3_fopinf.csv};
\addlegendentry{fi-OpInf}
\addplot [ultra thick, black, forget plot]
table {%
1 1e-06
1 1
};
\addplot [ultra thick, black, forget plot]
table {%
3 1e-06
3 1
};
\legend{};
\end{axis}

\end{tikzpicture}
	\end{subfigure}
	\hfill
	\begin{subfigure}[b]{0.3\textwidth}
		% This file was created with tikzplotlib v0.10.1.
\begin{tikzpicture}

\definecolor{darkgray176}{RGB}{176,176,176}
\definecolor{darkorange}{RGB}{255,140,0}
\definecolor{green01270}{RGB}{0,127,0}
\definecolor{lightgray204}{RGB}{204,204,204}

\begin{axis}[
width=\plotwtripple,
height=\plothtripple,
legend cell align={left},
legend style={
  fill opacity=0.8,
  draw opacity=1,
  text opacity=1,
  at={(0.97,0.03)},
  anchor=south east,
  draw=lightgray204
},
log basis y={10},
tick align=outside,
tick pos=left,
x grid style={darkgray176},
xlabel={$f$, Hz},
xmajorgrids,
xmin=-0.7395, xmax=15.7495,
xtick style={color=black},
y grid style={darkgray176},
ylabel={}, 
xmax=8, 
ymajorgrids,
ymin=1e-06,  
ymax=0.5,
yticklabel style = {font = \small, xshift = 7pt},
xticklabel style = {font = \small},
ymode=log,
ytick style={color=black},
ytick={1e-07,1e-06,1e-05,0.0001,0.001,0.01,0.1,1,10},
yticklabels={
  \(\displaystyle {10^{-7}}\),
  \(\displaystyle {10^{-6}}\),
  \(\displaystyle {10^{-5}}\),
  \(\displaystyle {10^{-4}}\),
  \(\displaystyle {10^{-3}}\),
  \(\displaystyle {10^{-2}}\),
  \(\displaystyle {10^{-1}}\),
  \(\displaystyle {10^{0}}\),
  \(\displaystyle {10^{1}}\)
}
]
\addplot [semithick, green01270, mark=triangle*, mark size=2.5, mark options={solid}, only marks]
table[x = fr, y = popinf,  col sep=comma] {pics/cbeam/csv/ex1_text_freq_3_5_popinf.csv};
\addlegendentry{p-OpInf}
\addplot [semithick, darkorange, mark=*, mark size=2.5, mark options={solid}, only marks]
table[x = fr, y = fopinf,  col sep=comma] {pics/cbeam/csv/ex1_text_freq_3_5_fopinf.csv};
\addlegendentry{fi-OpInf}
\addplot [ultra thick, black, forget plot]
table {%
3 1e-06
3 1
};
\addplot [ultra thick, black, forget plot]
table {%
5 1e-06
5 1
};
\legend{};
\end{axis}

\end{tikzpicture}
	\end{subfigure}
\vspace{-0.5cm}
	\caption{Cantilever beam example: ROMs of order $r = 4$, maximal relative error for different frequencies dependent on a frequency range test conditions $dt = 0.01s$, $T = 16s$.}
	\label{fig:cbeam_tend2_r4}
\end{figure}

The error behavior for the overhanging beam is shown in \Cref{fig:beam4_tend5_r3}.  We test the overhanging beam for $t \in [0 , 0.6]s$ with the time step $dt = 0.001s$ at different frequencies. The eigenfrequencies for this example are higher, correspondingly the reduced-order modeling and the tests were performed at higher frequencies.

From the plots, we can see that approximating the system behavior at higher frequencies is more challenging, but can be accomplished by using the snapshots for reduced-order modeling that correspond to the higher-frequency chirp signal. At the same time, the system behavior at frequencies lower than the lowest "learned" frequency is very likely to be correctly represented. Although both examples show that training at high frequencies guarantees the correct approximation of the system response at lower frequencies, it is not sufficient to generalize this observation. However, if the load conditions of a system are known, a good first guess for the training snapshots would be the data set corresponding to the highest frequency of interest.

\begin{figure}[h!]
	\centering
	\begin{subfigure}[b]{0.3\textwidth}
		% This file was created with tikzplotlib v0.10.1.
\begin{tikzpicture}

\definecolor{darkgray176}{RGB}{176,176,176}
\definecolor{darkorange}{RGB}{255,140,0}
\definecolor{green01270}{RGB}{0,127,0}
\definecolor{lightgray204}{RGB}{204,204,204}

\begin{axis}[
width=\plotwtripple,
height=\plothtripple,
legend cell align={left},
legend style={
  fill opacity=0.8,
  draw opacity=1,
  text opacity=1,
  at={(0.75,0.33)},
  anchor=south east,
  draw=lightgray204
},
log basis y={10},
tick align=outside,
tick pos=left,
x grid style={darkgray176},
xlabel={$f$, Hz},
xmajorgrids,
xmin=10, xmax=320,
xtick style={color=black},
y grid style={darkgray176},
ylabel={Max. relative error},
ymajorgrids,
ymin=1e-06, ymax=0.1,
ymode=log,
yticklabel style = {font = \small, xshift = 7pt},
xticklabel style = {font = \small},
ytick style={color=black},
ytick={1e-07,1e-06,1e-05,0.0001,0.001,0.01,0.1,1,10},
yticklabels={
  \(\displaystyle {10^{-7}}\),
  \(\displaystyle {10^{-6}}\),
  \(\displaystyle {10^{-5}}\),
  \(\displaystyle {10^{-4}}\),
  \(\displaystyle {10^{-3}}\),
  \(\displaystyle {10^{-2}}\),
  \(\displaystyle {10^{-1}}\),
  \(\displaystyle {10^{0}}\),
  \(\displaystyle {10^{1}}\)
}
]
\addplot [semithick, green01270, mark=triangle*, mark size=2.5, mark options={solid}, only marks]
table[x = fr, y = popinf,  col sep=comma] {pics/beam4/csv/ex2_text_freq_50_70_popinf.csv};
\addlegendentry{\texttt{p-OpInf}}
\addplot [semithick, darkorange, mark=*, mark size=2.5, mark options={solid}, only marks]
table[x = fr, y = fopinf,  col sep=comma] {pics/beam4/csv/ex2_text_freq_50_70_fopinf.csv};
\addlegendentry{\texttt{fi-OpInf}}
\addplot [ultra thick, black, forget plot]
table {%
50 1e-06
50 0.1
};
\addplot [ultra thick, black, forget plot]
table {%
70 1e-06
70 0.1
};
\end{axis}

\end{tikzpicture}
	\end{subfigure}
	\hspace{0.8cm}
	\begin{subfigure}[b]{0.3\textwidth}
		% This file was created with tikzplotlib v0.10.1.
\begin{tikzpicture}

\definecolor{darkgray176}{RGB}{176,176,176}
\definecolor{darkorange}{RGB}{255,140,0}
\definecolor{green01270}{RGB}{0,127,0}
\definecolor{lightgray204}{RGB}{204,204,204}

\begin{axis}[
width=\plotwtripple,
height=\plothtripple,
%legend cell align={left},
%legend style={
%  fill opacity=0.8,
%  draw opacity=1,
%  text opacity=1,
%  at={(0.97,0.03)},
%  anchor=south east,
%  draw=lightgray204
%},
log basis y={10},
tick align=outside,
tick pos=left,
x grid style={darkgray176},
xlabel={$f$, Hz},
xmajorgrids,
xmin=10, xmax=320,
xtick style={color=black},
y grid style={darkgray176},
ymajorgrids,
ymin=1e-06, ymax=0.1,
yticklabel style = {font = \small, xshift = 7pt},
xticklabel style = {font = \small},
ymode=log,
ytick style={color=black},
ytick={1e-07,1e-06,1e-05,0.0001,0.001,0.01,0.1,1,10},
yticklabels={
  \(\displaystyle {10^{-7}}\),
  \(\displaystyle {10^{-6}}\),
  \(\displaystyle {10^{-5}}\),
  \(\displaystyle {10^{-4}}\),
  \(\displaystyle {10^{-3}}\),
  \(\displaystyle {10^{-2}}\),
  \(\displaystyle {10^{-1}}\),
  \(\displaystyle {10^{0}}\),
  \(\displaystyle {10^{1}}\)
}
]
\addplot [semithick, green01270, mark=triangle*, mark size=2.5, mark options={solid}, only marks]
table[x = fr, y = popinf,  col sep=comma] {pics/beam4/csv/ex2_text_freq_140_160_popinf.csv};
%\addlegendentry{p-OpInf}
\addplot [semithick, darkorange, mark=*, mark size=2.5, mark options={solid}, only marks]
table[x = fr, y = fopinf,  col sep=comma] {pics/beam4/csv/ex2_text_freq_140_160_fopinf.csv};
%\addlegendentry{fi-OpInf}
\addplot [ultra thick, black, forget plot]
table {%
140 1e-06
140 0.1
};
\addplot [ultra thick, black, forget plot]
table {%
160 1e-06
160 0.1
};
\end{axis}

\end{tikzpicture}
	\end{subfigure}
	\hfill
	\begin{subfigure}[b]{0.3\textwidth}
		% This file was created with tikzplotlib v0.10.1.
\begin{tikzpicture}

\definecolor{darkgray176}{RGB}{176,176,176}
\definecolor{darkorange}{RGB}{255,140,0}
\definecolor{green01270}{RGB}{0,127,0}
\definecolor{lightgray204}{RGB}{204,204,204}

\begin{axis}[
width=\plotwtripple,
height=\plothtripple,
legend cell align={left},
legend style={
  fill opacity=0.8,
  draw opacity=1,
  text opacity=1,
  at={(0.97,0.03)},
  anchor=south east,
  draw=lightgray204
},
log basis y={10},
tick align=outside,
tick pos=left,
x grid style={darkgray176},
xlabel={$f$, Hz},
xmajorgrids,
xmin=10, xmax=320,
xtick style={color=black},
y grid style={darkgray176},
ymajorgrids,
ymin=1e-06, ymax=0.1,
ymode=log,
yticklabel style = {font = \small, xshift = 7pt},
xticklabel style = {font = \small},
ytick style={color=black},
ytick={1e-07,1e-06,1e-05,0.0001,0.001,0.01,0.1,1,10},
yticklabels={
  \(\displaystyle {10^{-7}}\),
  \(\displaystyle {10^{-6}}\),
  \(\displaystyle {10^{-5}}\),
  \(\displaystyle {10^{-4}}\),
  \(\displaystyle {10^{-3}}\),
  \(\displaystyle {10^{-2}}\),
  \(\displaystyle {10^{-1}}\),
  \(\displaystyle {10^{0}}\),
  \(\displaystyle {10^{1}}\)
}
]
\addplot [semithick, green01270, mark=triangle*, mark size=2.5, mark options={solid}, only marks]
table[x = fr, y = popinf,  col sep=comma] {pics/beam4/csv/ex2_text_freq_250_270_popinf.csv};
%\addlegendentry{p-OpInf}
\addplot [semithick, darkorange, mark=*, mark size=2.5, mark options={solid}, only marks]
table[x = fr, y = fopinf,  col sep=comma] {pics/beam4/csv/ex2_text_freq_250_270_fopinf.csv};
%\addlegendentry{fi-OpInf}
\addplot [ultra thick, black, forget plot]
table {%
250 1e-06
250 0.1
};
\addplot [ultra thick, black, forget plot]
table {%
270 1e-06
270 0.1
};
\end{axis}

\end{tikzpicture}
	\end{subfigure}
	\vspace{-0.7cm}
	\caption{Overhanging beam example: ROM of order $r = 3$, maximal relative error for different frequencies dependent on a frequency range test conditions $dt = 0.0002s$, $T = 1s$.}
	\label{fig:beam4_tend5_r3}
\end{figure}
\FloatBarrier

\newpage 
\subsection{Rotor}

\begin{figure}[h!] 
	\centering
	\includegraphics[scale = 0.4]{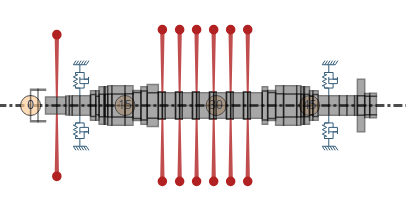}
	\caption{Rotor finite element model of dimension $n = 224$~\cite{ross2020}.}
	\label{fig:ex3}
\end{figure}

The next example is a rotor model, built and simulated using the Python library for rotordynamic analysis ROSS \cite{ross2020}. The model corresponds to the system \eqref{eq:original}, where the rotation speed $\Omega \neq 0$ and $\bG \neq 0$, i.e., the reduced-order modeling involves the identification of the gyroscopic term $\bG$. The rotor contains of 50 shaft elements, 7 disk elements, and 2 bearing elements, see \Cref{fig:ex3}. 
 \FloatBarrier
The total number of DOFs for this model is 224. Each node of the rotor model has 4 DOFs: translational displacements in two  directions ($\mathrm{x}_1$, $\mathrm{x}_2$) and two rotational displacements $\theta_1$ and $\theta_2$, i.e., $\textbf{x}^{\text{node}} = [\mathrm{x}_1, \; \mathrm{x}_2, \; \theta_1, \; \theta_2]$.
% about the $x_1$ and $x_2$ axes, corresponding to $\mathrm{x}_1$ and $\mathrm{x}_2$ respectively:
To generate the snapshots for the reduced-order modeling, we excite the rotor on one node in two directions with the chirp signal \eqref{eq:chirp_signal} with $\phi_0 = 0$ for $\mathrm{x}_1$-direction and $\phi_0 = 90$ for the $\mathrm{x}_2$-direction. We use the time step size $dt = 0.001 \si{\second}$, the rotation speed $\Omega_{\text{train}} = 600$\si{rad/s}, and force amplitude $\text{A}_\text{m} = 10 \si{\kilo\newton}$.  The corresponding external force and singular value decay are depicted in \Cref{fig:rotor_force}, where the frequency defined in \eqref{eq:chirp} starts at $f_0 = 1 \si{\hertz}$ and reaches $f_1 = 4 \si{\hertz}$ at the end-time point $T = 10\si{\second}$. 

%\begin{figure}[h!]
%	\begin{subfigure}{0.49\textwidth}
%		\input{pics/rotor/rotor_F_tend10_dt0p001_speed600_y.tex}
%		\caption{Input signal $\mathrm{u_m}(t)$ according to \eqref{eq:chirp} with a frequency range $[1, 4]$\si{\hertz}.}
%	\end{subfigure}
%\hfill
%	\begin{subfigure}{0.49\textwidth}
%		\input{pics/rotor/rotor_svd_tend10_dt0p001_speed600.tex}
%		\caption{SVD decay of the snapshot matrix $\bX$ that corresponds to the integration of the system \eqref{eq:original} using the input signal $\mathrm{u_m}(t)$ with a frequency range $[1, 4]$\si{\hertz}. The vertical line indicates the reduced order chosen for the model.}
%	\end{subfigure}
%	\label{fig:rotor_force}
%\end{figure}\FloatBarrier 

\begin{figure}[h!]
	\subcaptionbox{Input signal $\mathrm{u_m}(t)$ according to \eqref{eq:chirp} with a frequency range $[1, 4]$\si{\hertz}.}%
	[.45\linewidth]{\input{pics/rotor/rotor_F_tend10_dt0p001_speed600_y.tex}}
	\hfill
	\subcaptionbox{Singular value decay of the snapshot matrix $\bX$ that corresponds to the integration of the system \eqref{eq:original} using the input signal $\mathrm{u_m}(t)$ with a frequency range $[1, 4]$\si{\hertz}. The vertical line indicates the reduced order $r = 2$ chosen for the ROM.}
	[.49\linewidth]{% This file was created with tikzplotlib v0.10.1.
\begin{tikzpicture}

\definecolor{darkgray176}{RGB}{176,176,176}
\definecolor{steelblue31119180}{RGB}{31,119,180}
\definecolor{darkred}{RGB}{139,0,0}

\begin{axis}[
width=\plotw,
height=\ploth,
log basis y={10},
tick align=outside,
tick pos=left,
x grid style={darkgray176},
xlabel={r},
xmin=-4.95, xmax=83.95,
xtick style={color=black},
ymajorgrids,
xmajorgrids,
y grid style={darkgray176},
ylabel={\(\displaystyle \sigma / \sigma_{max}\)},
ymin=1.41253754462275e-18, ymax=7.07945784384138,
ymode=log,
ytick style={color=black},
ytick={1e-21,1e-18,1e-15,1e-12,1e-09,1e-06,0.001,1,1000,1000000},
yticklabels={
  \(\displaystyle {10^{-21}}\),
  \(\displaystyle {10^{-18}}\),
  \(\displaystyle {10^{-15}}\),
  \(\displaystyle {10^{-12}}\),
  \(\displaystyle {10^{-9}}\),
  \(\displaystyle {10^{-6}}\),
  \(\displaystyle {10^{-3}}\),
  \(\displaystyle {10^{0}}\),
  \(\displaystyle {10^{3}}\),
  \(\displaystyle {10^{6}}\)
},
xticklabel style={font = \small},
yticklabel style={font = \small}
]
\addplot [semithick, black, mark=*, mark size=2, mark options={solid}, only marks]
table {%
0 1
1 0.993510328178851
2 7.31954051266125e-05
3 6.8793078506506e-05
4 3.42728066017609e-05
5 9.21020101407251e-06
6 5.94207444229884e-07
7 6.19678606795975e-08
8 1.12066600720713e-08
9 9.43839301364663e-09
10 4.78215142277293e-09
11 1.02502799363527e-09
12 2.46071927034436e-10
13 1.43538277805822e-10
14 6.17901209868528e-11
15 2.27846846079977e-11
16 6.95875503268172e-12
17 2.79774307538069e-12
18 9.97698031622287e-13
19 3.95322474802994e-13
20 2.73542690246833e-13
21 1.15047303127722e-13
22 6.97235073103372e-14
23 1.43362784741599e-14
24 7.2602196662826e-15
25 4.62307973933235e-15
26 2.53565079548623e-15
27 1.46833058507905e-15
28 5.83041095872829e-16
29 4.82146387047161e-16
30 3.91324356702344e-16
31 1.74301242735614e-16
32 1.59619380021504e-16
33 9.95519841716537e-17
34 9.95519841716537e-17
35 9.95519841716537e-17
36 9.95519841716537e-17
37 9.95519841716537e-17
38 9.95519841716537e-17
39 9.95519841716537e-17
40 9.95519841716537e-17
41 9.95519841716537e-17
42 9.95519841716537e-17
43 9.95519841716537e-17
44 9.95519841716537e-17
45 9.95519841716537e-17
46 9.95519841716537e-17
47 9.95519841716537e-17
48 9.95519841716537e-17
49 9.95519841716537e-17
50 9.95519841716537e-17
51 9.95519841716537e-17
52 9.95519841716537e-17
53 9.95519841716537e-17
54 9.95519841716537e-17
55 9.95519841716537e-17
56 9.95519841716537e-17
57 9.95519841716537e-17
58 9.95519841716537e-17
59 9.95519841716537e-17
60 9.95519841716537e-17
61 9.95519841716537e-17
62 9.95519841716537e-17
63 9.95519841716537e-17
64 9.95519841716537e-17
65 9.95519841716537e-17
66 9.95519841716537e-17
67 9.95519841716537e-17
68 9.95519841716537e-17
69 9.95519841716537e-17
70 9.95519841716537e-17
71 9.95519841716537e-17
72 9.95519841716537e-17
73 9.95519841716537e-17
74 9.95519841716537e-17
75 9.95519841716537e-17
76 9.95519841716537e-17
77 9.95519841716537e-17
78 9.95519841716537e-17
79 9.95519841716537e-17
80 9.95519841716537e-17
81 9.95519841716537e-17
82 9.95519841716537e-17
83 9.95519841716537e-17
84 9.95519841716537e-17
85 9.95519841716537e-17
86 9.95519841716537e-17
87 9.95519841716537e-17
88 9.95519841716537e-17
89 9.95519841716537e-17
90 9.95519841716537e-17
91 9.95519841716537e-17
92 9.95519841716537e-17
93 9.95519841716537e-17
94 9.95519841716537e-17
95 9.95519841716537e-17
96 9.95519841716537e-17
97 9.95519841716537e-17
98 9.95519841716537e-17
99 9.95519841716537e-17
};
\addplot [ultra thick, darkred, dash pattern=on 14.8pt off 6.4pt]
table {%
2 1e-17
2 1
};
\end{axis}

\end{tikzpicture} \vspace{0.2cm}}
	\caption{External load and output snapshot characteristics of the rotor model.}
		\label{fig:rotor_force}
\end{figure}
\begin{figure}[h!]
	\centering
	\begin{subfigure}[b]{0.3\textwidth}
		%red
	\tikzexternalenable%
	\tikzsetnextfilename{pics/rotor/rotor_error_freq1_2_dt_0p001_tend10_r2_valid}%
	\filemodCmp{pics/rotor/rotor_error_freq1_2_dt_0p001_tend10_r2_valid.tex}{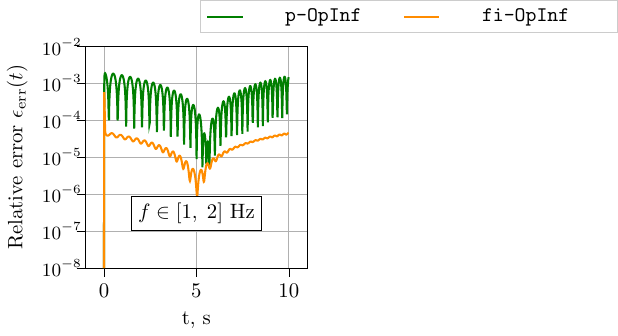}%
	{\tikzset{external/remake next}}{}%
	% This file was created with tikzplotlib v0.10.1.
\begin{tikzpicture}

\definecolor{darkgray176}{RGB}{176,176,176}
\definecolor{darkorange}{RGB}{255,140,0}
\definecolor{green01270}{RGB}{0,127,0}
\definecolor{lightgray204}{RGB}{204,204,204}

\begin{axis}[
width=\plotwtripple,
height=\plothtripple,
legend cell align={left},
legend style={
	fill opacity=0.8,
	draw opacity=1,
	text opacity=1,
	at={(2.4,1.06)},
	anchor=south east,
	draw=lightgray204,
	legend columns=2,
	nodes={inner xsep=20pt} 
},
log basis y={10},
tick align=outside,
tick pos=left,
x grid style={darkgray176},
xlabel={t, s},
xmajorgrids,
xtick style={color=black},
y grid style={darkgray176},
ylabel={Relative error $\epsilon_{\mathrm{err}}(t)$},
ymajorgrids,
ymin=1e-08, ymax=0.01,
ymode=log,
ytick style={color=black},
yticklabel style={color=black, font=\small, xshift = 6pt},
ytick={1e-08,1e-07,1e-06,1e-05,0.0001,0.001,0.01,0.1,1},
yticklabels={
  \(\displaystyle {10^{-8}}\),
  \(\displaystyle {10^{-7}}\),
  \(\displaystyle {10^{-6}}\),
  \(\displaystyle {10^{-5}}\),
  \(\displaystyle {10^{-4}}\),
  \(\displaystyle {10^{-3}}\),
  \(\displaystyle {10^{-2}}\),
  \(\displaystyle {10^{-1}}\),
  \(\displaystyle {10^{0}}\)
}
]
\addplot [line width=\linew,each nth point={15}, green01270]
table[x = time, y expr={\thisrow{popinf} <= \pgfkeysvalueof{/pgfplots/ymin} ? \pgfkeysvalueof{/pgfplots/ymin} : \thisrow{popinf}},  col sep=comma] {pics/rotor/csv/ex3_valid_error_freq_1_2.csv};
\addlegendentry{\texttt{p-OpInf}}
\addplot [line width=\linew, darkorange]
table[x = time, y expr={\thisrow{fopinf} <= \pgfkeysvalueof{/pgfplots/ymin} ? \pgfkeysvalueof{/pgfplots/ymin} : \thisrow{fopinf}}, each nth point={15}, col sep=comma] {pics/rotor/csv/ex3_valid_error_freq_1_2.csv};
\addlegendentry{\texttt{fi-OpInf}}
\node at (axis cs:5,0.0000001)  [draw, rectangle, fill = white, anchor=south] {$f \in [1, \; 2]$ Hz};
\end{axis}

\end{tikzpicture}%
	\tikzexternaldisable%

	\end{subfigure}
	\hspace{0.7cm}
	\begin{subfigure}[b]{0.3\textwidth}
		%red
	\tikzexternalenable%
	\tikzsetnextfilename{pics/rotor/rotor_error_freq2_3_dt_0p001_tend10_r2_valid}%
	\filemodCmp{pics/rotor/rotor_error_freq2_3_dt_0p001_tend10_r2_valid.tex}{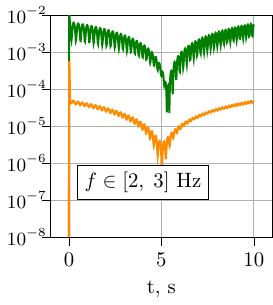}%
	{\tikzset{external/remake next}}{}%
	% This file was created with tikzplotlib v0.10.1.
\begin{tikzpicture}

\definecolor{darkgray176}{RGB}{176,176,176}
\definecolor{darkorange}{RGB}{255,140,0}
\definecolor{green01270}{RGB}{0,127,0}
\definecolor{lightgray204}{RGB}{204,204,204}

\begin{axis}[
width=\plotwtripple,
height=\plothtripple,
%legend cell align={left},
%legend style={
%  fill opacity=0.8,
%  draw opacity=1,
%  text opacity=1,
%  at={(0.97,0.03)},
%  anchor=south east,
%  draw=lightgray204
%},
log basis y={10},
tick align=outside,
tick pos=left,
x grid style={darkgray176},
xlabel={t, s},
xmajorgrids,
xtick style={color=black},
y grid style={darkgray176},
ymajorgrids,
ymin=1e-08, ymax=0.01,
ymode=log,
ytick style={color=black},
yticklabel style={color=black, font=\small, xshift = 6pt},
ytick={1e-08,1e-07,1e-06,1e-05,0.0001,0.001,0.01,0.1,1},
yticklabels={
  \(\displaystyle {10^{-8}}\),
  \(\displaystyle {10^{-7}}\),
  \(\displaystyle {10^{-6}}\),
  \(\displaystyle {10^{-5}}\),
  \(\displaystyle {10^{-4}}\),
  \(\displaystyle {10^{-3}}\),
  \(\displaystyle {10^{-2}}\),
  \(\displaystyle {10^{-1}}\),
  \(\displaystyle {10^{0}}\)
}
]
\addplot [line width=\linew,each nth point={15}, green01270]
table[x = time,y expr={\thisrow{popinf} <= \pgfkeysvalueof{/pgfplots/ymin} ? \pgfkeysvalueof{/pgfplots/ymin} : \thisrow{popinf}},  col sep=comma] {pics/rotor/csv/ex3_valid_error_freq_2_3.csv};
%\addlegendentry{p-OpInf}
\addplot [line width=\linew, darkorange]
table[x = time, y expr={\thisrow{fopinf} <= \pgfkeysvalueof{/pgfplots/ymin} ? \pgfkeysvalueof{/pgfplots/ymin} : \thisrow{fopinf}}, each nth point={15},  col sep=comma] {pics/rotor/csv/ex3_valid_error_freq_2_3.csv};
%\addlegendentry{fi-OpInf}
\node at (axis cs:4,0.0000001)  [draw, rectangle, fill = white, anchor=south] {$f \in [2, \; 3]$ Hz};
\end{axis}

\end{tikzpicture}%
	\tikzexternaldisable%

	\end{subfigure}
	\hfill
	\begin{subfigure}[b]{0.3\textwidth}
		%red
	\tikzexternalenable%
	\tikzsetnextfilename{pics/rotor/rotor_error_freq1_4_dt_0p001_tend10_r2_valid}%
	\filemodCmp{pics/rotor/rotor_error_freq1_4_dt_0p001_tend10_r2_valid.tex}{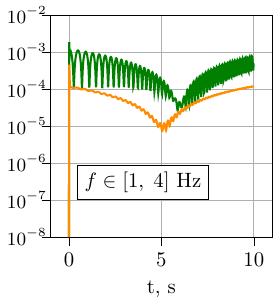}%
	{\tikzset{external/remake next}}{}%
	% This file was created with tikzplotlib v0.10.1.
\begin{tikzpicture}

\definecolor{darkgray176}{RGB}{176,176,176}
\definecolor{darkorange}{RGB}{255,140,0}
\definecolor{green01270}{RGB}{0,127,0}
\definecolor{lightgray204}{RGB}{204,204,204}

\begin{axis}[
width=\plotwtripple,
height=\plothtripple,
%legend cell align={left},
%legend style={
%  fill opacity=0.8,
%  draw opacity=1,
%  text opacity=1,
%  at={(0.97,0.03)},
%  anchor=south east,
%  draw=lightgray204
%},
log basis y={10},
tick align=outside,
tick pos=left,
x grid style={darkgray176},
xlabel={t, s},
xmajorgrids,
xtick style={color=black},
y grid style={darkgray176},
ymajorgrids,
ymin=1e-08, ymax=0.01,
ymode=log,
ytick style={color=black},
yticklabel style={color=black, font=\small, xshift = 6pt},
ytick={1e-08,1e-07,1e-06,1e-05,0.0001,0.001,0.01,0.1,1},
yticklabels={
  \(\displaystyle {10^{-8}}\),
  \(\displaystyle {10^{-7}}\),
  \(\displaystyle {10^{-6}}\),
  \(\displaystyle {10^{-5}}\),
  \(\displaystyle {10^{-4}}\),
  \(\displaystyle {10^{-3}}\),
  \(\displaystyle {10^{-2}}\),
  \(\displaystyle {10^{-1}}\),
  \(\displaystyle {10^{0}}\)
}
]
\addplot [line width=\linew,each nth point={15}, green01270]
table[x = time, y expr={\thisrow{popinf} <= \pgfkeysvalueof{/pgfplots/ymin} ? \pgfkeysvalueof{/pgfplots/ymin} : \thisrow{popinf}},  col sep=comma] {pics/rotor/csv/ex3_valid_error_freq_1_4.csv};
%\addlegendentry{p-OpInf}
\addplot [line width=\linew,each nth point={15}, darkorange]
table[x = time, y expr={\thisrow{fopinf} <= \pgfkeysvalueof{/pgfplots/ymin} ? \pgfkeysvalueof{/pgfplots/ymin} : \thisrow{fopinf}},  col sep=comma] {pics/rotor/csv/ex3_valid_error_freq_1_4.csv};
%\addlegendentry{fi-OpInf}
\node at (axis cs:4,0.0000001)  [draw, rectangle, fill = white, anchor=south] {$f \in [1, \; 4]$ Hz};
\end{axis}

\end{tikzpicture}%
	\tikzexternaldisable%

	\end{subfigure}
	\vspace{-0.7cm}
	\caption{Rotor example: Relative error vs. time for the ROMs of dimension $r = 2$.}
	\label{fig:rotor_valid_errors}
\end{figure}

 We also investigated different frequency ranges for snapshot generation. \Cref{fig:rotor_valid_errors} plots the relative errors over time for the resulting ROMs, both of which show reasonable accuracy throughout the simulation. Varying the training frequency has a negligible impact on \texttt{fi-OpInf}, which outperforms the parameterized operator inference across all considered snapshot sets.
 \FloatBarrier
 Unlike previous structural beam models, single-frequency simulations show no significant change in maximum error for the force-informed operator inference model. On the other hand, the parameterized operator inference approach, is highly sensitive to the training frequency range and shows large errors for frequencies larger than 1 $\si{\hertz}$ even inside the "training window", see \Cref{fig:rotor_test_errors_a}, \Cref{fig:rotor_test_errors_b}. Broadening the training frequency range initially improved the maximum relative errors, as shown in \Cref{fig:rotor_test_errors_c}, but further extending the range to 8 $\si{\hertz}$ yielded no improvements at higher frequencies. This behavior of the parametrized operator inference method can likely be attributed to the increased number of parameters and a high sensitivity to initial guesses.
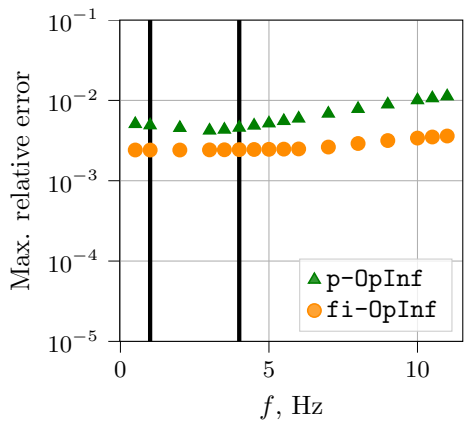
\begin{figure}[h!]
	%\hspace{-1.7cm}
	\begin{subfigure}[b]{0.3\textwidth}
		% This file was created with tikzplotlib v0.10.1.
\begin{tikzpicture}

\definecolor{darkgray176}{RGB}{176,176,176}
\definecolor{darkorange}{RGB}{255,140,0}
\definecolor{green01270}{RGB}{0,127,0}
\definecolor{lightgray204}{RGB}{204,204,204}

\begin{axis}[
width=\plotwtripple,
height=\plothtripple,
legend cell align={left},
legend style={
	fill opacity=0.8,
	draw opacity=1,
	text opacity=1,
	at={(0.97,0.03)},
	anchor=south east,
	draw=lightgray204
},
log basis y={10},
tick align=outside,
tick pos=left,
x grid style={darkgray176},
xlabel={$f$, Hz},
xmajorgrids,
xmin=-0.025, xmax=8,
xtick style={color=black},
y grid style={darkgray176},
ylabel={Max. relative error},
ymajorgrids,
ymin=1e-05, ymax=0.1,
ymode=log,
ytick style={color=black},
ytick={1e-07,1e-06,1e-05,0.0001,0.001,0.01,0.1,1,10},
yticklabels={
  \(\displaystyle {10^{-7}}\),
  \(\displaystyle {10^{-6}}\),
  \(\displaystyle {10^{-5}}\),
  \(\displaystyle {10^{-4}}\),
  \(\displaystyle {10^{-3}}\),
  \(\displaystyle {10^{-2}}\),
  \(\displaystyle {10^{-1}}\),
  \(\displaystyle {10^{0}}\),
  \(\displaystyle {10^{1}}\)
},
xticklabel style={font = \small},
yticklabel style={font = \small}
]
\addplot [semithick, green01270, mark=triangle*, mark size=2.5, mark options={solid}, only marks]
table[x = fr, y = popinf,  col sep=comma] {pics/rotor/csv/ex3_text_freq_1_2.csv};
\addlegendentry{\texttt{p-OpInf}}
\addplot [semithick, darkorange, mark=*, mark size=2.5, mark options={solid}, only marks]
table[x = fr, y = fopinf,  col sep=comma] {pics/rotor/csv/ex3_text_freq_1_2.csv};
\addlegendentry{\texttt{fi-OpInf}}
\addplot [ultra thick, black, forget plot]
table {%
1 1e-05
1 0.1
};
\addplot [ultra thick, black, forget plot]
table {%
2 1e-05
2 0.1
};
\end{axis}

\end{tikzpicture}
		\label{fig:rotor_test_errors_a}
	\end{subfigure}
	\hspace{0.8cm}
	%\hfill
	\begin{subfigure}[b]{0.3\textwidth}
		% This file was created with tikzplotlib v0.10.1.
\begin{tikzpicture}

\definecolor{darkgray176}{RGB}{176,176,176}
\definecolor{darkorange}{RGB}{255,140,0}
\definecolor{green01270}{RGB}{0,127,0}
\definecolor{lightgray204}{RGB}{204,204,204}

\begin{axis}[
width=\plotwtripple,
height=\plothtripple,
legend cell align={left},
legend style={
	fill opacity=0.8,
	draw opacity=1,
	text opacity=1,
	at={(0.97,0.03)},
	anchor=south east,
	draw=lightgray204
},
log basis y={10},
tick align=outside,
tick pos=left,
x grid style={darkgray176},
xlabel={$f$, Hz},
xmajorgrids,
xmin=-0.025, xmax=8,
xtick style={color=black},
y grid style={darkgray176},
ymajorgrids,
ymin=1e-05, ymax=0.1,
ymode=log,
ytick style={color=black},
ytick={1e-07,1e-06,1e-05,0.0001,0.001,0.01,0.1,1,10},
yticklabels={
  \(\displaystyle {10^{-7}}\),
  \(\displaystyle {10^{-6}}\),
  \(\displaystyle {10^{-5}}\),
  \(\displaystyle {10^{-4}}\),
  \(\displaystyle {10^{-3}}\),
  \(\displaystyle {10^{-2}}\),
  \(\displaystyle {10^{-1}}\),
  \(\displaystyle {10^{0}}\),
  \(\displaystyle {10^{1}}\)
},
xticklabel style={font = \small},
yticklabel style={font = \small}
]
\addplot [semithick, green01270, mark=triangle*, mark size=2.5, mark options={solid}, only marks]
table[x = fr, y = popinf,  col sep=comma] {pics/rotor/csv/ex3_text_freq_2_3.csv};
\addlegendentry{\texttt{p-OpInf}}
\addplot [semithick, darkorange, mark=*, mark size=2.5, mark options={solid}, only marks]
table[x = fr, y = fopinf,  col sep=comma] {pics/rotor/csv/ex3_text_freq_2_3.csv};
\addlegendentry{\texttt{fi-OpInf}}
\addplot [ultra thick, black, forget plot]
table {%
2 1e-05
2 0.1
};
\addplot [ultra thick, black, forget plot]
table {%
3 1e-05
3 0.1
};
\end{axis}

\end{tikzpicture}
		\label{fig:rotor_test_errors_b}
	\end{subfigure}
	\hfill
	\begin{subfigure}[b]{0.3\textwidth}
		% This file was created with tikzplotlib v0.10.1.
\begin{tikzpicture}

\definecolor{darkgray176}{RGB}{176,176,176}
\definecolor{darkorange}{RGB}{255,140,0}
\definecolor{green01270}{RGB}{0,127,0}
\definecolor{lightgray204}{RGB}{204,204,204}

\begin{axis}[
width=\plotwtripple,
height=\plothtripple,
legend cell align={left},
legend style={
	fill opacity=0.8,
	draw opacity=1,
	text opacity=1,
	at={(0.97,0.03)},
	anchor=south east,
	draw=lightgray204
},
log basis y={10},
tick align=outside,
tick pos=left,
x grid style={darkgray176},
xlabel={$f$, Hz},
xmajorgrids,
xmin=-0.025, xmax=8,
xtick style={color=black},
y grid style={darkgray176},
ymajorgrids,
ymin=1e-05, ymax=0.1,
ymode=log,
ytick style={color=black},
ytick={1e-07,1e-06,1e-05,0.0001,0.001,0.01,0.1,1,10},
yticklabels={
  \(\displaystyle {10^{-7}}\),
  \(\displaystyle {10^{-6}}\),
  \(\displaystyle {10^{-5}}\),
  \(\displaystyle {10^{-4}}\),
  \(\displaystyle {10^{-3}}\),
  \(\displaystyle {10^{-2}}\),
  \(\displaystyle {10^{-1}}\),
  \(\displaystyle {10^{0}}\),
  \(\displaystyle {10^{1}}\)
},
xticklabel style={font = \small},
yticklabel style={font = \small}
]
\addplot [semithick, green01270, mark=triangle*, mark size=2.5, mark options={solid}, only marks]
table[x = fr, y = popinf,  col sep=comma] {pics/rotor/csv/ex3_text_freq_1_4.csv};
\addlegendentry{\texttt{p-OpInf}}
\addplot [semithick, darkorange, mark=*, mark size=2.5, mark options={solid}, only marks]
table[x = fr, y = fopinf,  col sep=comma] {pics/rotor/csv/ex3_text_freq_1_4.csv};
\addlegendentry{\texttt{fi-OpInf}}
\addplot [ultra thick, black, forget plot]
table {%
1 1e-05
1 0.1
};
\addplot [ultra thick, black, forget plot]
table {%
4 1e-05
4 0.1
};
\end{axis}

\end{tikzpicture}
		\label{fig:rotor_test_errors_c}
	\end{subfigure}
	\vspace{-0.7cm}
	\caption{Rotor example: Maximal relative error for different monofrequent simulations of the ROM of order $r = 2$ constructed using the snapshots at the different frequency range. Test simulations performed with $dt = 0.001s$, end time $T = 5s$.}
	\label{fig:rotor_test_errors}
\end{figure} \FloatBarrier The ROM constructed from the snapshots that correspond to the chirp input signal with the frequency range $[1,4] \si{\hertz}$ is able to accurately reproduce the system trajectories at different frequencies outside and inside the training range. We compare the spatial 3D and 2D trajectories in \Cref{fig:rotor_test_errors} that show a good agreement with the FOM results.
\begin{figure}[h!]
	\hspace{-1cm}
	\begin{subfigure}[b]{0.45\textwidth}
		%\includesvg[scale=0.5]{rotor3d.pdf}
		\includesvg[scale=0.7]{./pics/rotor/rotor_3D} 
		\vspace{0.3cm}
		\caption{3D trajectory.}	
	\end{subfigure}
	%	\hspace{1.1cm}
	\hfill
	\begin{subfigure}[b]{0.45\textwidth}
		% This file was created with tikzplotlib v0.10.1.
\begin{tikzpicture}

\definecolor{darkgray176}{RGB}{176,176,176}
\definecolor{darkorange}{RGB}{255,140,0}
\definecolor{green}{RGB}{0,128,0}
\definecolor{lightgray204}{RGB}{204,204,204}
\definecolor{steelblue31119180}{RGB}{31,119,180}

\begin{axis}[
width=\plotw,
height=\ploth,
legend cell align={left},
legend style={
  fill opacity=0.8,
  draw opacity=1,
  text opacity=1,
  at={(0.91,0.5)},
  anchor=east,
  draw=lightgray204
},
tick align=outside,
tick pos=left,
x grid style={darkgray176},
xlabel={Displacement $\mathrm{x}_1$, mm},
xmajorgrids,
%xmin=-8.57771360896477e-08, xmax=8.57944486163233e-08,
xtick style={color=black},
%xtick={-1e-07,-7.5e-08,-5e-08,-2.5e-08,0,2.5e-08,5e-08,7.5e-08,1e-07},
%xticklabels={
%  \ensuremath{-}10.0,
%  \ensuremath{-}7.5,
%  \ensuremath{-}5.0,
%  \ensuremath{-}2.5,
%  0.0,
%  2.5,
%  5.0,
%  7.5,
%  10.0
%},
y grid style={darkgray176},
ylabel={Displacement $\mathrm{x}_2$, mm},
ymajorgrids,
%ymin=-8.61390717171374e-08, ymax=8.58179563599021e-08,
ytick style={color=black},
%ytick={-1e-07,-7.5e-08,-5e-08,-2.5e-08,0,2.5e-08,5e-08,7.5e-08,1e-07},
%yticklabels={
%  \ensuremath{-}10.0,
%  \ensuremath{-}7.5,
%  \ensuremath{-}5.0,
%  \ensuremath{-}2.5,
%  0.0,
%  2.5,
%  5.0,
%  7.5,
%  10.0
%},
every y tick scale label/.style={
	at={(-0.07,0.99)},xshift=12pt,anchor=north west,inner sep=0pt
},
every x tick scale label/.style={
	at={(0.9,0.0)},xshift=12pt,anchor=north west,inner sep=0pt
},
xticklabel style={font = \small},
yticklabel style={font = \small}
]
\addplot [semithick, steelblue31119180] table[x = fomX, y = fomY,  col sep=comma] {pics/rotor/csv/ex3_2D_traj_freq_1_4.csv};
\addlegendentry{\texttt{FOM}}
\addplot [semithick, darkorange, dash pattern=on 9pt off 30pt on 3pt off 30pt] table[x = fopinfX, y = fopinfY,  col sep=comma] {pics/rotor/csv/ex3_2D_traj_freq_1_4.csv};
\addlegendentry{\texttt{fi-OpInf}}
\addplot [semithick, green, dash pattern=on 3pt off 30pt] table[x = popinfX, y = popinfY,  col sep=comma] {pics/rotor/csv/ex3_2D_traj_freq_1_4.csv};
\addlegendentry{\texttt{p-OpInf}}
\end{axis}

\end{tikzpicture}
		\vspace{-0.5cm}
		\caption{2D trajectory.}
		
	\end{subfigure}
	%\vspace{-0.5cm}
	\caption{Rotor example: Comparison of the spatial trajectories of FOM and ROM of order $r = 2$ constructed using the chirp input signal with the frequency range $[1, 4] \si{\hertz}$. Test simulations performed at the frequency $3 \si{\hertz}$ with $dt = 0.001s$, end time $T = 20s$.}
	\label{fig:rotor_3Dand2D}
\end{figure} \FloatBarrier
 
For rotor dynamics, analyzing the reduced-order model under different rotational speeds $\Omega$ is of primary interest. In practice, the underlying system matrices are dependent on the rotation speed, therefore a parameter-dependent inference approach is required. The derivation of this parameterized model is left for future research.

In summary, the numerical results for the given three examples show a good consistency of the \texttt{p-OpInf}, \texttt{fi-OpInf} models and the original model simulations. Constructing the ROM using only one simulation data set that corresponds to a harmonic frequency-variable input signal already leads to trajectories with reasonable accuracy for a wide range of other harmonic input frequencies. In most cases, the accuracy of \texttt{fi-OpInf} outperforms \texttt{p-OpInf} ROM. One of the possible reasons is the choice of the initial values for the unknown coefficients during the \texttt{p-OpInf} optimization.

\section{Conclusions}

In this work, we presented a new methodology (\texttt{p-OpInf}) to identify the reduced operators of a typical second-order in time ODE system from input and state data. The key step of our approach is the solution of a nonlinear optimization problem, which can be efficiently implemented and solved using automatic differentiation tools such as PyTorch. 
Furthermore, we complement the methodology with a parametrization of the unknown reduced system operators to guarantee the preservation of the mathematical properties of the system matrices, i.e., the symmetric positive (semi-)definiteness of the mass, stiffness and damping matrices, as well as the skew-symmetric properties of the gyroscopic term for rotating structures.
We demonstrated the numerical performance of \texttt{p-OpInf} for 3 examples. Two structural beam models are used to test the \texttt{p-OpInf} surrogate model for structural vibrations at different frequencies, and the rotor model shows how the \texttt{p-OpInf} reduced model approximates the rotational behavior by identifying the gyroscopic system matrix. We compared the results with the existing force-informed operator inference methodology (\texttt{fi-OpInf}), which has also been extended to include rotational structures.
The numerical experiments show acceptable accuracy for a very small reduced-order and a wide range of loading conditions, using only one set of simulation data for the reduced-order model construction. However, achieving a good approximation of the high-frequency behavior is challenging. In this respect, the \texttt{fi-OpInf} method shows much more reliable error behavior than the \texttt{p-OpInf} method, i.e., lower error for more loading cases. On the other hand, \texttt{fi-OpInf} modeling requires knowledge of the force positions and distributions, which are assumed to be unknown and can be identified by the \texttt{p-OpInf} method.
The error behavior of \texttt{p-OpInf} can be improved by using better initial guesses for the system operators. However, this aspect is not investigated in this work. Another open question is also the choice of optimal simulation conditions for the data set used during the reduced-order modeling. \newline

\textbf{Data availability statement} \\
The simulation data, and code to produce the results
presented in this paper are openly available at \url{https://doi.org/10.5281/zenodo.18064656}. The numerical
experiments were implemented in Python 3.11.5, and run on a Lenovo Laptop with a 11th Gen Intel®Core™ i7-11800H 32 GB RAM. The system was running Linux Mint  version 22 (64Bit). 
%The procedure described in Section is implemented in Python 3.11.5, using PyTorch
\newline

\textbf{Acknowledgments} \\
The authors acknowledge the support and computational resources provided by the Max Planck Institute for Dynamics of Complex Technical Systems, Magdeburg. This research was also supported by the Research Training Group “Mathematical Complexity Reduction”, which is a Graduiertenkolleg (DFG-GRK 2297) funded by Deutsche Forschungsgemeinschaft (DFG), Germany.
The authors thank Vahid Nateghi and Dr. Jens Saak for reviewing the code.

%\section{Acknowledgments}
%The authors acknowledge the support and computational resources provided by the Max Planck Institute for Dynamics of Complex Technical Systems, Magdeburg. This research was also supported by a Research Training Group - Mathematical Complexity Reduction, which is a Graduiertenkolleg (DFG-GRK 2297) funded by Deutsche Forschungsgemeinschaft (DFG).

\newpage
\addcontentsline{toc}{section}{References}
\bibliographystyle{plainurl}

\bibliography{popinf_lit}

%\newpage
%\appendix
%\input{real_numres.tex}
  
\end{document}